\documentclass[11pt]{article}
\usepackage[title]{appendix}

\usepackage{amsmath,amsfonts,amssymb,algorithm,algorithmic,psfrag}
\usepackage{amsthm}
\usepackage{array}
\usepackage{psfrag}
\usepackage{subfigure}
\usepackage{bm}
\usepackage{color}
\usepackage[curve]{xypic}
\usepackage{bigstrut}

\newtheorem{theorem}{Theorem}[section]
\newtheorem{lemma}[theorem]{Lemma}
\newtheorem{corollary}[theorem]{Corollary}

\theoremstyle{definition}
\newtheorem{definition}[theorem]{Definition}

\theoremstyle{remark}

\numberwithin{equation}{section}

\usepackage{xypic}

\usepackage{color} 
\usepackage{ifpdf}
\ifpdf
    \usepackage[pdftex]{graphicx}
    \usepackage[pdftex]{hyperref}
    \hypersetup{
        unicode=false,          
        pdftoolbar=true,        
        pdfmenubar=true,        
        pdffitwindow=false,     
        pdfstartview={FitH},    
        pdftitle={AKG Article},      
        pdfauthor={Andrew Gillette},   
        pdfsubject={Mathematics},    
        pdfcreator={Andrew Gillette},  
        pdfproducer={Andrew Gillette}, 
        pdfkeywords={math, mathematics}, 
        pdfnewwindow=true,      
        colorlinks=true,        
        linkcolor=red,          
        citecolor=blue,         
        filecolor=magenta,      
        urlcolor=cyan           
    }

\else
    \usepackage{graphicx}
    \usepackage{pstricks}
    
    \newcommand{\href}[2]{#2}
\fi

\usepackage[mathscr]{euscript}
\usepackage{amsbsy}

\usepackage[shortlabels]{enumitem}
\usepackage{textcomp} 
\usepackage[strict]{changepage} 
\usepackage{xypic}

\usepackage{amssymb}

\DeclareMathOperator{\myspan}{span}

\DeclareMathOperator{\tr}{tr}

\newcommand{\bbN}{\mathbb{N}}

\newcommand{\bbR}{\mathbb{R}}

\newcommand{\calH}{\mathcal{H}}

\newcommand{\calJ}{\mathcal{J}}

\newcommand{\calP}{\mathcal{P}}
\newcommand{\calQ}{\mathcal{Q}}

\newcommand{\calS}{\mathcal{S}}

\newcommand{\beq}{\begin{equation}}
\newcommand{\eeq}{\end{equation}}

\newcommand{\ds}{\displaystyle}

\newcommand{\hdiv}{H(\textnormal{div})}
\newcommand{\hcurl}{H(\textnormal{curl})}

\newcommand{\ul}[1]{\underline{#1}}
\newcommand{\p}{\partial}
\newcommand{\raw}{\rightarrow}

\newcommand{\ldeg}{\textnormal{ldeg}}

\newcommand{\nedelec}{N{\'e}d{\'e}lec}

\newcommand{\nfwdtab}[2]{
#1 & #2
}

\newcommand{\nfrdtab}[3]{
#1 & #2 & #3
}

\theoremstyle{plain}

\renewcommand{\emptyset}{\varnothing}

\title{Computational Serendipity and Tensor Product \\ 
Finite Element Differential Forms}

\author{Andrew Gillette\footnote{
Department of Mathematics,
University of Arizona,
Tucson, AZ 85721. \hfill
{\it agillette@math.arizona.edu}}
~~~~Tyler Kloefkorn\footnote{
AAAS Science \& Technology Policy Fellow, hosted at NSF, Alexandria, VA. \hfill
{\it tyler.kloefkorn@aaas-fpi.com}}
~~~~Victoria Sanders\footnote{
Department of Mathematics,
University of Arizona,
Tucson, AZ 85721. \hfill
{\it victoriasanders@email.arizona.edu}}
}
\date{}

\begin{document}
\maketitle

\begin{abstract}
Many conforming finite elements on squares and cubes are elegantly classified into families by the language of finite element exterior calculus and presented in the Periodic Table of the Finite Elements. 
Use of these elements varies, based principally on the ease or difficulty in finding a ``computational basis'' of shape functions for element families. The tensor product family, $\calQ^-_r\Lambda^k$, is most commonly used because computational basis functions are easy to state and implement.
The trimmed and non-trimmed serendipity families, $\calS^-_r\Lambda^k$ and $\calS_r\Lambda^k$ respectively, are used less frequently because they are newer to the community and, until now, lacked a straightforward technique for computational basis construction.
This represents a missed opportunity for computational efficiency as the serendipity elements in general have fewer degrees of freedom than elements of equivalent accuracy from the tensor product family. 
Accordingly, in pursuit of easy adoption of the serendipity families, we present complete lists of computational bases for both serendipity families, for any order $r\geq 1$ and for any differential form order $0\leq k\leq n$, for problems in dimension $n=2$ or $3$.
The bases are defined via shared subspace structures, allowing easy comparison of elements across families.
We use and include code in \textit{SageMath} to find, list, and verify these  computational basis functions.
\end{abstract}

\maketitle

\section{Introduction}
\label{sec:intro}

Three major families of conforming finite element differential forms on cubical meshes have been classified using the framework of finite element exterior calculus: $\calQ_r^-\Lambda^k$, $\calS_r\Lambda^k$, and $\calS_r^-\Lambda^k$, called the tensor product, serendipity, and trimmed serendipity families, respectively.
The tensor product family includes the well-known \nedelec~elements for $\hcurl$ and $\hdiv$ problems and has been used to great success in, for instance, computational electromagnetism.
The widespread use of the tensor product family is due, in large part, to the straightforward computational basis functions available for its implementation. 
The two kinds of serendipity families have gained attention in recent years for their prospect of attaining equivalent rates of convergence as their tensor product counterpart but at significantly reduced computational expense. 
Nevertheless, these benefits have not been realized in practice due primarily to the lack of suitable computational basis functions required for their implementation.  

We present here, for the first time, a comprehensive list of basis functions for all three families of finite element differential forms, in a unified notational framework for dimensions 2 and 3, of any form order $k$ and of any polynomial order $r\geq 1$. 
Our list provides a clear association of basis functions to vertices, edges, faces, and the interior of squares and cubes, following a formal definition of computational basis functions.
We highlight the similarities and differences among the constructions, especially how certain `building blocks' of spaces of functions are used in multiple families.  
Identifying these blocks aids in efficient coding while also illuminating the hierarchical structure of the spaces.

To understand the subtleties required for compiling such a list, consider the construction of a computational basis for $\calS_1\Lambda^1(\square_3)$, a serendipity space of differential 1-forms on the cube $\square_3:=[-1,1]^3$.
The dimension of the space is 24, with 2 degrees of freedom associated to each edge of the cube.
The degrees of freedom require that a basis for $\calS_1\Lambda^1(\square_3)$ provide constant and linear order approximation on each edge of the cube. 
Hence, it is practical to seek basis functions with canonical associations to edges.
Functions providing constant order approximation are straightforward:\begin{equation}
\label{eq:obv}
\left\{(y\pm 1)(z\pm 1)dx,~~(x\pm 1)(z\pm 1)dy,~~(x\pm 1)(y\pm 1)dz\right\},
\end{equation}
a total of 12 functions.
Observe that each function corresponds to an edge of the cube in a clear fashion.
For instance, the function $(y+1)(z+1)dx$, when restricted to an edge of $\square_3$, is zero on every edge except $\{y=1, z=1\}$ where it has value 4, thereby providing approximation of constants on that edge.\footnote{In traditional vector calculus notation, $(y+1)(z+1)dx$ can be thought of as the 3-vector $\begin{bmatrix} (y+1)(z+1) & 0 & 0 \end{bmatrix}^T$.}
To attain linear order approximation on each edge, an obvious guess to complete the basis would be the additional 12 functions
\begin{equation}
\label{eq:obv-gues}
\left\{x(y\pm 1)(z\pm 1)dx,~~y(x\pm 1)(z\pm 1)dy,~~z(x\pm 1)(y\pm 1)dz\right\}.
\end{equation}
Unfortunately, none of the forms from (\ref{eq:obv-gues}) are in $\calS_1\Lambda^1(\square_3)$! More importantly, the span of (\ref{eq:obv}) together with (\ref{eq:obv-gues}) is \ul{not} the space $\calS_1\Lambda^1(\square_3)$.
One issue with the above guess lies in the fact that the highest order coefficient -- $xyz$ in each case -- cannot appear in isolation, e.g.\ $xyz\ dx\not\in\calS_1\Lambda^1(\square_3)$.
To get $xyz$ as a coefficient of, say, $dx$, it is necessary to also include some polynomial coefficients in $dy$ and $dz$.
For instance, using the definition of $\calS_1\Lambda^1(\square_3)$, we can check directly that
\begin{equation}
\label{eq:1fex-prep}
2xyz\ dx + x^2z\ dy + x^2y\ dz \in \calS_1\Lambda^1(\square_3).
\end{equation}
Similarly, we can show that the forms
\[\begin{array}{rcrcr}
2xy\ dx &+& x^2\ dy &+& 0\ dz, \\
2xz\ dx &+& 0 \ dy  &+& x^2\ dz, \\
2x\ dx &+& (-z -1)dy &+& (-y-1) dz
\end{array}\]
are also elements of $\calS_1\Lambda^1(\square_3)$.
Taking the sum of these three forms together with the form appearing in (\ref{eq:1fex-prep}) and factoring yields 
\begin{equation}
\label{eq:1fex}
2x(y+1)(z+1)\ dx + (z+1)(x^2-1)\ dy + (y+1)(x^2 - 1)\ dz,
\end{equation}
which has the desired polynomial in front of $dx$ (with coefficient 2), along with terms in $dy$ and $dz$ that have value 0 when restricted to any edge of the cube.
Thus, we take the form from (\ref{eq:1fex}) and its 11 analogous counterparts (3 other sign changes for the $(\ast+1)$ terms and 8 other forms with variables permuted) to complete a basis for $\calS_1\Lambda^1(\square_3)$ with canonical associations as desired.

The example of $\calS_1\Lambda^1(\square_3)$ is illustrative of the general challenge in constructing computational bases for serendipity style families on squares and cubes.
For an order $r$ element, the ``obvious'' functions are suitable up to at least order $r-1$ approximation.  
For order $r$ approximation, the ``obvious'' computational functions may not be in the approximation space on their own, but each one differs from a function in the space only by a function that vanishes identically on the boundary of the domain (i.e. a bubble function) or by a function that provides a similar kind of local approximation. 
Interestingly, we find that the serendipity and trimmed serendipity spaces exhibit opposite behavior in regards to which form types require the addition of bubble functions to complete a computational basis.

The remainder of the paper is organized as follows.
We begin in Section~\ref{sec:bkgd} with a review of relevant background and notation from finite element exterior calculus, including the definitions of regular and trimmed serendipity spaces.
In Section~\ref{sec:bases-theory}, we define computational bases for the $\calQ_r^-\Lambda^k$, $\calS_r\Lambda^k$ and $\calS^-_r\Lambda^k$ families on squares and cubes, in unified notation, and discuss some of the patterns that emerge from this perspective.
Next, in Section~\ref{sec:bases-sage}, we describe code we have written using the DifferentialForms package in \textit{SageMath} to generate the computational bases for any of the three families, given only the inputs $n$ (2 or 3), $r\geq 1$, and $0\leq k\leq n$.
We also use this code to verify that the bases we define are, indeed, bases for the desired spaces using a simple linear algebra technique.
A portion of the code is available at \url{http://math.arizona.edu/~agillette/research/computationalBases/}.
A discussion of applications and future directions for this work is given in Section~\ref{sec:conc}.
Appendix~\ref{app:lists} gives complete lists of some of the serendipity and trimmed serendipity bases for $n=3$.

\section{Notation and Relation to Prior Work}
\label{sec:bkgd}

We briefly review relevant notation from the finite element exterior calculus literature~\cite{AFW2006,AFW2010}, primarily following conventions used for the the tensor product family~\cite{ABB2012}, the serendipity family~\cite{AA2014}, the trimmed serendipity family~\cite{GK2016}, and the Periodic Table of the Finite Elements~\cite{AL2014}.

\subsection*{Spaces of forms on $\bbR^n$.}
Define the \textit{form monomial} $x^\alpha dx_\sigma$ to be the differential $k$-form on $\bbR^n$ given by
\begin{equation}
\label{eq:form-mon}
x^{\alpha}dx_{\sigma}:=\left(x_1^{\alpha_1}x_2^{\alpha_2}\dots x_n^{\alpha_n}\right)dx_{\sigma(1)}\wedge\dots\wedge dx_{\sigma(k)},
\end{equation}
where $\alpha\in\bbN^n$ is a multi-index and $\sigma$ is a subset of $\{1,\ldots,n\}$ consisting of $k$ distinct elements $\sigma(1),\ldots,\sigma(k)$ with $0\leq k\leq n$.
The $x^\alpha$ portion of a form monomial is called its \textit{coefficient} and $dx_\sigma$ is called its \textit{alternator}. 
For convenience of notation, wedges will often be omitted from expressions involving forms.

For the tensor product spaces, we need a notion of maximum degree per variable.
The space
\begin{equation}
\label{eq:qrm-part}
\left[\bigotimes_{i=1}^n\calP_{r-\delta_{i,\sigma}}\right]dx_{\sigma(1)}\wedge\dots\wedge dx_{\sigma(k)}
\end{equation}
refers to the span of form monomials on $\bbR^n$ whose monomial coefficient is at most degree $r$ in variables that do not appear in its alternator and at most degree $r-1$ in variables that do appear in its alternator.
Note that $0$-forms have an empty alternator and thus if $k=0$, (\ref{eq:qrm-part}) is simply the standard scalar-valued tensor product space of order $r$.

For non-tensor product spaces, we need other measures of degree.
The \textit{degree} of $x^{\alpha}dx_{\sigma}$ is $|\alpha|:=\sum_{i=1}^n\alpha_i$.
Let $\sigma^\ast$ denote the complement of $\sigma$, i.e. $\sigma^\ast:=\{1,\ldots,n\}-\sigma$.
The \textit{linear degree} of $x^{\alpha}dx_{\sigma}$ is defined to be
\begin{equation}
\label{eq:lindeg-def}
\ldeg(x^\alpha dx_\sigma):=\#\{i\in\sigma^\ast~:~\alpha_i=1\}.
\end{equation}
Put differently, the linear degree of $x^\alpha dx_\sigma$ counts the number of entries in $\alpha$ equal to 1, excluding entries whose indices appear in $\sigma$.
Note that if $k=0$ then $\sigma=\varnothing$ and there is no `exclusion' in the counting of linear degree.
Likewise, if $k=n$ then $\sigma^\ast=\varnothing$ and $\ldeg(x^\alpha dx_\sigma)=0$ for any $\alpha$.
The linear degree of the sum of two or more form monomials is defined as the minimum of the linear degrees of the summands.

The space of differential $k$-forms with polynomial coefficients of homogeneous degree $r$ is denoted $\calH_r\Lambda^k(\bbR^n)$.
A basis for this space is the set of form monomials such that $|\alpha|=r$ and $|\sigma|=k$.
The exterior derivative $d$ and Koszul operator $\kappa$ are maps
\[d:\calH_r\Lambda^k(\bbR^n)\raw\calH_{r-1}\Lambda^{k+1}(\bbR^n)\quad \kappa:\calH_r\Lambda^k(\bbR^n)\raw\calH_{r+1}\Lambda^{k-1}(\bbR^n).\]
In coordinates, they are defined on form monomials by
\begin{align}
d(x^\alpha dx_\sigma) & := \sum_{i=1}^n \left(\frac{\partial x^\alpha}{\partial x_i}dx_i\right)\wedge dx_{\sigma(1)}\wedge\dots\wedge dx_{\sigma(k)}, \label{eq:def-d}\\
\kappa(x^\alpha dx_\sigma) & := \sum_{i=1}^k \left((-1)^{i+1}x^\alpha x_{\sigma(i)}\right)dx_{\sigma(1)}\wedge\cdots\wedge\widehat{dx_{\sigma(i)}}\wedge\cdots\wedge dx_{\sigma(k)}.\label{eq:def-kappa}
\end{align}
The notation $\widehat{dx_{\sigma(i)}}$ indicates that the term is omitted from the wedge product.

The subset of $\calH_r\Lambda^k(\bbR^n)$ that has linear degree at least $\ell$ is denoted
\begin{equation}
\label{eq:hrlk-def}
		\calH_{r,l}\Lambda^k(\bbR^n):=\left\{\omega\in\calH_r\Lambda^k(\bbR^n)~|~\mbox{ldeg }\omega\geq l\right\}.
\end{equation}
A key building block for both the serendipity and trimmed serendipity spaces is
\begin{equation}
\label{eq:jrk-def}
\calJ_r \Lambda^k(\bbR^n) := \sum_{\ell \geq 1} \kappa \calH_{r+\ell-1, \ell} \Lambda^{k+1}(\bbR^n).
\end{equation}

We can now define all the major families of finite elements in compatible terminology.
We omit including $(\bbR^n)$ throughout for ease of notation.
The spaces of polynomial and trimmed polynomial differential $k$-forms are then
\begin{align}
\calP_r\Lambda^k & :=\bigoplus_{j=0}^r\calH_j\Lambda^k, \label{eq:prk-def}\\
\calP_r^-\Lambda^k & := \calP_{r-1}\Lambda^k \oplus \kappa\calH_{r-1}\Lambda^{k+1}.
\end{align}
The tensor product, serendipity, and trimmed serendipity spaces are
\begin{align}
\calQ_r^-\Lambda^k & := \bigoplus_{\sigma\in\Sigma(k,n)}\left[\bigotimes_{i=1}^n\calP_{r-\delta_{i,\sigma}}\right]dx_{\sigma(1)}\wedge\dots\wedge dx_{\sigma(k)}
,\\
\calS_r\Lambda^k & := \calP_r\Lambda^k \oplus \calJ_r \Lambda^k\oplus  d \calJ_{r+1}\Lambda^{k-1}, \\
\calS_r^-\Lambda^k & := \calP_r^-\Lambda^k \oplus \calJ_r\Lambda^k \oplus d \calJ_{r} \Lambda^{k-1}. \label{eq:srmk-def}
\end{align}
The direct summation in $\calQ_r^-\Lambda^k$ is taken over $\Sigma(k,n)$, which denotes the space of increasing maps from $\{1,\ldots, k\}$ into $\{1,\ldots, n\}$.

Many of these spaces are the approximation spaces for elements that have been around for decades, including $\calQ_r^-\Lambda^1(\bbR^2)$ (Raviart-Thomas~\cite{RT1977}), $\calQ_r^-\Lambda^1(\bbR^3)$ and $\calQ_r^-\Lambda^2(\bbR^3)$ (\nedelec~\cite{N1980,N1986}), and $\calS_r\Lambda^1(\square_2)$  (Brezzi-Douglas-Marini~\cite{BDM85}).

\subsection*{Form spaces on cubical geometries.}
We now need some terminology to describe spaces of forms and their bases on cubical geometries.
Define the $n$-dimensional cube $\square_n:=[-1,1]^n$ for any $n\geq 1$.

The \textit{trace} of a differential $k$-form on a co-dimension 1 hyperplane $f\subset\bbR^n$ is the pullback of the form via the inclusion map $f\hookrightarrow\bbR^n$.
Let $x^\alpha dx_\sigma$ be a form monomial and let $f$ be the hyperplane defined by $x_i=c$ for some fixed $1\leq i\leq n$ and  constant $c$.
Then the trace of the form monomial on $f$ is
\[\tr_f (x^\alpha dx_\sigma)=\begin{cases} 0, & i\in\sigma,\\ \left(x^\alpha|_{x_i=c}\right)dx_\sigma, & i\not\in\sigma.
\end{cases}\]
Note that the trace function distributes over sums of form monomials.
Further, the composition of trace functions can be used to compute the trace of a form on hyperplanes of codimension greater than 1.
It follows from the definition that the trace of a $k$-form on a region of dimension $\ell<k$ is always 0.

Given a polynomial differential $k$-form $\omega$ and a cubical geometry $\square_n$, define
\begin{equation}
\label{eq:def-m-omega}
 m_\omega := \min_{f\prec\square_n}\left\{\dim f ~:~\tr_f\omega\not = 0\right\},
\end{equation}
where the min is taken over sub-faces of $\square_n$ of any dimension $0,\ldots n$.
The value of $m_\omega$ is interpreted as the minimum dimension of a face on which $\omega$ has non-zero trace.
For example, consider the 1-forms $\eta=(x+1)dy$ and $\zeta=(x^2-1)dy$ in the case $n=2$.
Both $\eta$ and $\zeta$ have zero trace on all vertices of $\square_2$ since they are 1-forms.
Observe that $\eta$ has trace $2\ dy$ on the edge $\{x=1\}$ and thus has $m_\eta=1$, while $\zeta$ has zero trace on every edge of the square, but is non-zero on the interior of the square, making $m_\zeta=2$. 

We can now define precisely what we mean by a computational basis.

\begin{definition}
A set of differential forms $\{\omega_i\}$ on a cubical geometry $\square_n$ is a \textit{computational basis} for a particular finite element differential form space $\mathcal X\Lambda^k(\square_n)$ if it satisfies the following conditions:
\begin{enumerate}
\item The set $\{\omega_i\}$ is a basis for $\mathcal X\Lambda^k(\square_n)$.
\item Each $\omega\in\{\omega_i\}$ is associated to a unique face $f\prec\square_n$ with $\dim f= m_\omega$ such that $\tr_g\omega=0$ for any face $g\prec\square_n$, $g\not=f$ with $\dim g=\dim f$.
\end{enumerate}
\end{definition}
A computational basis is a local basis in that it is defined on a single reference element.
Global basis functions can then be defined by stitching together local basis functions according to mesh connectivity.
In this way, a computational basis provides a \textit{geometric decomposition} of a differential form approximation space in that the basis elements are in a fixed  correspondence to the degrees of freedom.

Computational bases were provided for the spaces $\calP_r\Lambda^k$ and $\calP_r^-\Lambda^k$ on triangles, tetrahedra, and general $n$-simplices by Arnold, Falk, and Winther in~\cite{AFW2006gd}.
Their procedure for $\calP_r^-\Lambda^k(\triangle_n)$ is as follows: construct barycentric coordinates $\{\lambda_i\}$ associated to each vertex of the simplex $\triangle_n$, build Whitney $k$-forms from these functions, e.g.\ $\{\lambda_i d\lambda_j - \lambda_j d\lambda_i\}$ for $1$-forms, multiply the Whitney forms by $r-1$ additional functions from $\{\lambda_i\}$, and finally reduce the resulting spanning set to a basis.
Their procedure for $\calP_r\Lambda^k(\triangle_n)$ is similar.

While it is possible to extend this elegant technique from simplicial to cubical geometries, as attempted by Bossavit~\cite{B2010}, the procedure quickly becomes more involved and subtle due to the large number of possible vertex $k$-tuples in non-simplices. 
The work of Gillette, Rand and Bajaj~\cite{GRB2014} and Chen and Wang~\cite{CW2016} pursue this approach in the more general context of polygons and polyhedra.
Since we deal exclusively with squares and cubes in this work, we find it is more straightforward to define the basis functions directly rather than as linear combinations from a larger spanning set.

\section{Computational bases for $\calQ_r^-\Lambda^k$, $\calS_r\Lambda^k$, and $\calS_r^-\Lambda^k$}
\label{sec:bases-theory}

We state the bases for the three spaces of finite element differential forms in Table~\ref{tab:dim2} (for squares) and Table~\ref{tab:dim3} (for cubes).
The basis elements are grouped by the parameter $m$, which refers to the definition given in \eqref{eq:def-m-omega}.
Given a value of $0\leq k\leq n$, the basis for $\calQ_r^-\Lambda^k$, $\calS_r\Lambda^k$, or $\calS_r^-\Lambda^k$ is constructed by looking at the relevant portion of the table and taking the union of the spaces in the appropriate column of the table.
For instance,
\begin{equation}
\label{eq:S312-comp-bas-ex}
\calS_3\Lambda^1(\square_2) = \left(\bigoplus_{i=0}^{2} E_i\Lambda^1(\square_2)\right) \oplus \tilde{E}_3\Lambda^1(\square_2)\oplus\left(\bigoplus_{i=2}^{3} F_i\Lambda^1(\square_2)\right)
\end{equation}
The tables appear first, followed by definitions of all the subspaces that they reference.
We then give a discussion of patterns observed in the tables and strategies for implementation.

\begin{table}
\begin{tabular}{c|c|lll}
$n=2$& $m$ & $k = 0$ & $k = 1$ & $k = 2$ \\[1mm]
\hline
& 0 & $V\Lambda^0(\square_2)$ & $\emptyset$ & $\emptyset$ \\[2mm]
$\calQ^-_r\Lambda^k$ & 1 & $\ds\bigoplus_{i=0}^{r-2} E_i\Lambda^0(\square_2)$ & $\ds\bigoplus_{i=0}^{r-1} E_i\Lambda^1(\square_2)$ & $\emptyset$ \\[2mm]
& 2 & $\ds\bigoplus_{i=1}^{r-1} F^{\otimes}_i\Lambda^0(\square_2)$ & $\ds\bigoplus_{i=1}^{r-1} F^{\otimes}_i\Lambda^1(\square_2)$ & $\ds\bigoplus_{i=1}^r F^{\otimes}_i\Lambda^2(\square_2)$ \\[4mm]
\hline
& 0 & $V\Lambda^0(\square_2)$ & $\emptyset$ & $\emptyset$ \\[2mm]
$\calS_r\Lambda^k$ & 1 & $\ds\bigoplus_{i=0}^{r-2} E_i\Lambda^0(\square_2)$ & $\ds\bigoplus_{i=0}^{r-1} E_i\Lambda^1(\square_2) \oplus \tilde{E}_r\Lambda^1(\square_2)$ & $\emptyset$ \\[2mm]
& 2 & $\ds\bigoplus_{i=4}^{r} F_i\Lambda^0(\square_2)$ & $\ds\bigoplus_{i=2}^{r} F_i\Lambda^1(\square_2)$ & $\ds\bigoplus_{i=0}^r F_i\Lambda^2(\square_2)$ \\[4mm]
\hline
& 0 & $V\Lambda^0(\square_2)$ & $\emptyset$ & $\emptyset$ \\[2mm]
$\calS^-_r\Lambda^k$ & 1 & $\ds\bigoplus_{i=0}^{r-2} E_i\Lambda^0(\square_2)$ & $\ds\bigoplus_{i=0}^{r-1} E_i\Lambda^1(\square_2)$ & $\emptyset$ \\[2mm]
& 2 & $\ds\bigoplus_{i=4}^{r} F_i\Lambda^0(\square_2)$ & $\ds\bigoplus_{i=2}^{r-1} F_i\Lambda^1(\square_2) \oplus \tilde{F}_{r}\Lambda^1(\square_2)$ & $\ds\bigoplus_{i=0}^{r-1} F_i\Lambda^2(\square_2)$ \\[4mm]
\hline
\end{tabular}
\caption{Computational bases in two dimensions according to form order $k$.  The parameter $m$ indicates the dimension of the geometry to which the basis functions in a row are associated.}
\label{tab:dim2}
\end{table}

\begin{table}[h]
\begin{adjustwidth}{-.2in}{-.5in}
\begin{tabular}{c|c|llll}
$n=3$& $m$ & $k = 0$ & $k = 1$ & $k = 2$ & $k=3$ \\[1mm]
\hline
& 0 & $V\Lambda^0(\square_3)$ & $\emptyset$ & $\emptyset$ & $\emptyset$ \\[2mm]
$\calQ^-_r\Lambda^k$
& 1 & $\ds\bigoplus_{i=0}^{r-2} E_i\Lambda^0(\square_3)$ & $\ds\bigoplus_{i=0}^{r-1} E_i\Lambda^1(\square_3)$ & $\emptyset$ & $\emptyset$  \\[2mm]
& 2 & $\ds\bigoplus_{i=1}^{r-1} F^{\otimes}_i\Lambda^0(\square_3)$ & $\ds\bigoplus_{i=1}^{r-1} F^{\otimes}_i\Lambda^1(\square_3)$ & $\ds\bigoplus_{i=1}^{r} F^{\otimes}_i\Lambda^2(\square_3)$ & $\emptyset$ \\[2mm] 
& 3 & $\ds\bigoplus_{i=1}^{r-1} I^{\otimes}_i\Lambda^0(\square_3)$ & $\ds\bigoplus_{i=1}^{r-1} I^{\otimes}_i\Lambda^1(\square_3)$ & $\ds\bigoplus_{i=1}^{r-1} I^{\otimes}_i\Lambda^2(\square_3)$ & $\ds\bigoplus_{i=1}^{r-1} I^{\otimes}_i\Lambda^3(\square_3)$ \\[4mm]
\hline
& 0 &  $V\Lambda^0(\square_3)$ & $\emptyset$ & $\emptyset$ & $\emptyset$ \\[2mm]
$\calS_r\Lambda^k$ 
& 1 & $\ds\bigoplus_{i=0}^{r-2} E_i\Lambda^0(\square_3)$ & $\ds\bigoplus_{i=0}^{r-1} E_i\Lambda^1(\square_3) \oplus \tilde{E}_r\Lambda^1(\square_3)$ & $\emptyset$ & $\emptyset$ \\[2mm]
& 2 & $\ds\bigoplus_{i=4}^{r} F_i\Lambda^0(\square_3)$ & $\ds\bigoplus_{i=2}^{r-1} F_i\Lambda^1(\square_3) \oplus \hat{F}_{r}\Lambda^1(\square_3)$ & $\ds\bigoplus_{i=0}^{r-1} F_i\Lambda^2(\square_3) \oplus \tilde F_r\Lambda^2(\square_3)$ & $\emptyset$  \\[2mm] 
& 3 &  $\ds\bigoplus_{i=6}^{r} I_i\Lambda^0(\square_3)$ & $\ds\bigoplus_{i=4}^{r} I_i\Lambda^1(\square_3)$ &  $\ds\bigoplus_{i=2}^r I_i\Lambda^2(\square_3)$ & $\ds\bigoplus_{i=2}^r I_i\Lambda^3(\square_3)$ \\[4mm]
\hline
& 0 & $V\Lambda^0(\square_3)$ & $\emptyset$ & $\emptyset$  & $\emptyset$ \\[2mm]
$\calS^-_r\Lambda^k$ 
& 1 & $\ds\bigoplus_{i=0}^{r-2} E_i\Lambda^0(\square_3)$ & $\ds\bigoplus_{i=0}^{r-1} E_i\Lambda^1(\square_3)$ & $\emptyset$ & $\emptyset$ \\[2mm]
& 2 & $\ds\bigoplus_{i=4}^{r} F_i\Lambda^0(\square_3)$ & $\ds\bigoplus_{i=2}^{r-1} F_i\Lambda^1(\square_3) \oplus \tilde{F}_{r}\Lambda^1(\square_3)$ & $\ds\bigoplus_{i=0}^{r-1} F_i\Lambda^2(\square_3)$  & $\emptyset$ \\[2mm]
& 3 & $\ds\bigoplus_{i=6}^{r} I_i\Lambda^0(\square_3)$ & $\ds\bigoplus_{i=4}^{r-1} I_i\Lambda^1(\square_3)\oplus \tilde I_r\Lambda^1(\square_3)$ & $\ds\bigoplus_{i=2}^{r-1} I_i\Lambda^2(\square_3) \oplus \tilde{I}_r\Lambda^2(\square_3)$ & $\ds\bigoplus_{i=2}^{r-1} I_i\Lambda^3(\square_3)$ \\[4mm]
\hline
\end{tabular}
\text{}\\
\end{adjustwidth}
\caption{Computational bases in three dimensions according to form order $k$.  The parameter $m$ indicates the dimension of the geometry to which the basis functions in a row are associated.}
\label{tab:dim3}
\end{table}

\noindent
\subsection*{Two dimensions}
The 0-form spaces are constructed from the following spaces:
\[
\begin{array}{rll}
\text{space} & \text{functions} & \text{$i$ values} \\[2mm]
\hline\\[-2mm]
V\Lambda^0(\square_2) &  \left\{ (x \pm 1)(y \pm 1) \right\} & \\[2mm]
E_i\Lambda^0(\square_2) &  \left\{ y^i(x \pm 1)(y^2-1), ~~x^i(y \pm 1)(x^2 - 1) \right\} & i\geq 0 \\[2mm]
F_i\Lambda^0(\square_2) &  \left\{ x^jy^k(x^2-1)(y^2-1) ~:~ j+k = i-4 \right\} & i\geq 4\\[2mm]
F^{\otimes}_i\Lambda^0(\square_2) & \left\{ \quad\quad\quad\textnormal{---\quad\textquotesingle\textquotesingle\quad---}\quad\quad\quad  ~:~ \max(j, k) = i-1\} \right\} & i\geq 1 \\[2mm]
\end{array}
\]

\newpage
The 1-form spaces in 2D are constructed from the following spaces.
Each function in the sets listed below should be interpreted as the polynomial coefficient in the first column in front of $dx$ plus the polynomial in the second column in front of $dy$.
\[
\begin{array}{rll}
\text{space} & \text{functions listed by $\left\{\begin{array}{ll} dx & dy\end{array}\right\}$} & \text{$i$ values} \\[2mm]
\hline\\[-2mm]
E_i\Lambda^1(\square_2) & \left\{\begin{array}{ll} \nfwdtab {x^i(y \pm 1)}0 \\ \nfwdtab 0{y^i(x \pm 1)} \\
\end{array}\right\} & i\geq 0\\[8mm]
\tilde{E}_i\Lambda^1(\square_2) & 
\left\{\begin{array}{ll}  \nfwdtab {(i+1)x^i(y \pm 1)} {x^{i-1}(x^2 - 1)} \\ \nfwdtab {y^{i-1}(y^2 - 1)} {(i+1)y^i(x \pm 1)} \end{array}\right\} & i\geq 0\\[8mm]
F_i\Lambda^1(\square_2) & \left\{\begin{array}{ll} \nfwdtab {x^jy^k(y^2 - 1)}0 \\ \nfwdtab 0{x^jy^k(x^2 - 1)}\end{array}~:~j + k= i - 2\right\} & i\geq 2\\[8mm]
F^{\otimes}_i\Lambda^1((\square_2) & 
\left\{\quad\qquad\qquad\textnormal{---\quad\textquotesingle\textquotesingle\quad---}\quad\qquad\qquad ~:~ \max(j, k-1) = i - 1\right\} & i\geq 1\\[8mm]
\tilde{F}_i\Lambda^1(\square_2) & 
\left\{\begin{array}{lll} 
\nfwdtab {y^{i-2}(y^2 - 1)} 0 \\ 
\nfwdtab 0{x^{i-2}(x^2 - 1)} \\ 
\nfwdtab {x^ky^{i-k-2}(y^2-1)} {-x^{k-1}y^{i-k-1}(x^2-1)} &~:~ 1 \leq k \leq i-2\\
\end{array} \right\} & i\geq 2 
\end{array}
\]

The 2-form spaces in 2D are constructed from the following spaces.
\[
\begin{array}{rll}
\text{space} & \text{functions} & \text{$i$ values} \\[2mm]
\hline\\[-2mm]
F_i\Lambda^2(\square_2) & \left\{ x^jy^k\ dxdy ~:~ j + k = i \right\} & i\geq 0 \\[3mm]
F^{\otimes}_i\Lambda^2(\square_2) & \left\{ \quad\textnormal{---~\textquotesingle\textquotesingle~---}\quad ~:~ \max(j, k) = i - 1 \right\} & i\geq 1
\end{array}
\]

\subsection*{Three dimensions}
The 0-form spaces are constructed from the following spaces:
\[
\begin{array}{rll}
\text{space} & \text{functions} & \text{$i$ values} \\[2mm]
\hline\\[-2mm]
V\Lambda^0(\square_3) &  \left\{ (x \pm 1)(y \pm 1)(z\pm 1) \right\} & \\[2mm]
E_i\Lambda^0(\square_3) &  \left\{ \begin{array}{l} z^i(x\pm 1)(y\pm 1)(z^2-1) \\ y^i(x \pm 1)(z\pm 1)(y^2-1) \\ x^i(y \pm 1)(z \pm 1)(x^2 - 1) \end{array}\right\} & i\geq 0 \\[8mm]
F_i\Lambda^0(\square_3) &  \left\{ \begin{array}{c} x^jy^k(z\pm 1)(x^2-1)(y^2-1) \\ x^jz^k(y\pm 1)(x^2-1)(z^2-1) \\ y^jz^k(x\pm 1)(y^2-1)(z^2-1) \end{array} ~:~ j+k = i-4\right\} & i\geq 4\\[8mm]
F^{\otimes}_i\Lambda^0(\square_3) & \left\{ \quad\qquad\qquad\textnormal{---\quad\textquotesingle\textquotesingle\quad---}\qquad\qquad\quad  ~:~  \max(j, k) = i-1 \right\} & i\geq 1 \\[8mm]
I_i\Lambda^0(\square_3) &  \left\{ x^jy^kz^\ell(x^2-1)(y^2-1)(z^2-1) : j+k+\ell = i-6 \right\} & i\geq 6 \\[2mm]
I^{\otimes}_i\Lambda^0(\square_3) & \left\{ \quad\quad\qquad\textnormal{---\quad\textquotesingle\textquotesingle\quad---}\qquad\quad\quad  ~:~  \max(j,k,\ell) = i-1 \right\} & i\geq 1
\end{array}
\]

The 1-form spaces in 3D are constructed from the following spaces.
Each function in the sets listed below should be interpreted as the polynomial coefficient in the first column in front of $dx$ plus the polynomial in the second column in front of $dy$ plus the polynomial in the third column in front of $dz$.

\begin{adjustwidth}{-.8in}{-.9in}
$
\begin{array}{rl}
\begin{array}{r}
\text{space and} \\
\text{$i$ values}
\end{array}
& \text{functions listed by $\left\{\begin{array}{lll} dx & dy & dz \end{array}\right\}$}  \\[2mm]
 \\
\hline\\[-2mm]
\begin{array}{c}
E_i\Lambda^1(\square_3)\\
\hfill i\geq 0
\end{array} & \left\{\begin{array}{lll}
\nfrdtab {x^i(y \pm 1)(z \pm 1)}{0}{0} \\ \nfrdtab 0{y^i(x \pm 1)(z\pm 1)}0 \\ \nfrdtab 00{z^i(x\pm 1)(y\pm 1)}
\end{array}\right\} \\[8mm]
\begin{array}{c}\tilde{E}_i\Lambda^1(\square_3) \\
\hfill i\geq 0
\end{array} 
 & \left\{\begin{array}{lll}
\nfrdtab {(i+1)x^i(y \pm 1)(z \pm 1)}{x^{i-1}(z \pm 1)(x^2-1)}{x^{i-1}(y \pm 1)(x^2-1)} \\
\nfrdtab{y^{i-1}(z \pm 1)(y^2-1)}{(i+1)y^i(x \pm 1)(z \pm 1)}{y^{i-1}(x \pm 1)(y^2-1)} \\
\nfrdtab{z^{i-1}(y \pm 1)(z^2-1)}{z^{i-1}(x \pm 1)(z^2-1)}{(i+1)z^i(x \pm 1)(y \pm 1)}\\
\end{array}\right\} \\[6mm]
& \text{}\quad\text{In the above definition, there are exactly two independent sign options per line,} \\
& \text{}\quad\text{e.g.\ if $(y-1)$ is selected in the $dx$ term, then $(y-1)$ is selected for the $dz$ term.}\\
[4mm]
\begin{array}{c}
F_i\Lambda^1(\square_3) \\
\hfill i\geq 2
\end{array} 
& \left\{\begin{array}{lll}
\nfrdtab {x^jy^k(z\pm 1)(y^2-1)}{0}{0}\\
\nfrdtab {x^jz^k(y\pm 1)(z^2-1)}{0}{0}\\ 
\nfrdtab {0}{y^jx^k(z\pm 1)(x^2-1)}{0}\\  
\nfrdtab {0}{y^jz^k(x\pm 1)(z^2-1)}{0}\\ 
\nfrdtab {0}{0}{z^jy^k(y\pm 1)(x^2-1)}\\ 
\nfrdtab {0}{0}{z^jx^k(x\pm 1)(y^2-1)}
\end{array}~: ~j + k= i - 2 \right\}  \\[16mm]
\begin{array}{c}
F^{\otimes}_i\Lambda^1((\square_3) \\
\hfill i\geq 1
\end{array} 
& 
\left\{~\qquad\qquad\qquad\qquad\qquad\qquad\textnormal{---\quad\textquotesingle\textquotesingle\quad---}\qquad\qquad\qquad\qquad\qquad\qquad~  ~:~ \max(k, j-1) = i - 1\right\} \\[8mm]
\begin{array}{c}
\hat{F}_i\Lambda^1(\square_3) \\
\hfill i\geq 2
\end{array} 
& \left\{\begin{array}{lll}
\nfrdtab {y^{i-2}(z\pm 1)(y^2-1)}{0}{0}\\
\nfrdtab {z^{i-2}(y\pm 1)(z^2-1)}{0}{0}\\ 
\nfrdtab 0{x^{i-2}(z\pm 1)(x^2-1)}0 \\ 
\nfrdtab 0{z^{i-2}(x\pm 1)(z^2 - 1)}0 \\
\nfrdtab 00{x^{i-2}(y\pm 1)(x^2-1)}\\ 
\nfrdtab 00{y^{i-2}(x\pm 1)(y^2-1)} \\
\nfrdtab {(i+1)x^jy^{i-j-2}(z\pm 1)(y^2-1)}{0}{x^{j-1}y^{i-j-2}(x^2-1)(y^2-1)} \\
\nfrdtab {(i+1)x^jz^{i-j-2}(y\pm 1)(z^2-1)}{x^{j-1}z^{i-j-2}(x^2-1)(z^2-1)}{0}\\ 
\nfrdtab {y^{j-1}z^{i-j-2}(y^2-1)(z^2-1)}{(i+1)y^jz^{i-j-2}(x\pm 1)(z^2-1)}{0}\\
\nfrdtab {0}{(i+1)y^jx^{i-j-2}(z\pm 1)(x^2-1)}{y^{j-1}x^{i-j-2}(x^2-1)(y^2-1)} \\
\nfrdtab {0}{z^{j-1}x^{i-j-2}(x^2-1)(z^2-1)}{(i+1)z^jx^{i-j-2}(y\pm 1)(x^2-1)} \\
\nfrdtab {z^jy^{i-j-2}(y^2-1)(z^2-1)}{0}{(i+1)z^jy^{i-j-2}(x\pm 1)(y^2-1)} \\
\hline
\nfrdtab {}{~:~1\leq j\leq i-2}{}\\
\end{array}\right\}  \\[32mm]
\begin{array}{c}
\tilde{F}_i\Lambda^1(\square_3) \\
\hfill i\geq 2
\end{array}
& 
\left\{\begin{array}{llll}
\multicolumn{3}{c}{\text{first 6 rows are the same as those of $\hat{F}_i\Lambda^1(\square_3)$}} \\
\nfrdtab {x^jy^{i-j-2}(z\pm 1)(y^2-1)}{-x^{j-1}y^{i-j-1}(z\pm 1)(x^2-1)}0 \\
\nfrdtab {x^jz^{i-j-2}(y\pm 1)(z^2 - 1)}0{-x^{j-1}z^{i-j-1}(y\pm 1)(x^2-1)}\\
\nfrdtab 0{y^jz^{i-j-2}(x\pm 1)(z^2-1)}{-y^{j-1}z^{i-j-1}(x\pm 1)(y^2 -1)} \\
\hline
\nfrdtab {}{~:~1\leq j\leq i-2}{}\\
\end{array}
\right\} \\[16mm]
\end{array}
$
\end{adjustwidth}
\newpage
\begin{adjustwidth}{-.8in}{-.8in}
$
\begin{array}{rl}
\begin{array}{r}
\text{space and} \\
\text{$i$ values}
\end{array}
& \text{functions listed by $\left\{\begin{array}{lll} dx & dy & dz \end{array}\right\}$}  \\[2mm]
 \\
\hline\\[-2mm]
\begin{array}{c}
I_i\Lambda^1(\square_3)\\
\hfill i\geq 4
\end{array}
 & \left\{\begin{array}{lll}
 \nfrdtab {x^jy^kz^\ell(y^2- 1)(z^2-1)}00 \\  
 \nfrdtab 0{x^jy^kz^\ell(x^2-1)(z^2-1)}0 \\ 
 \nfrdtab 00{x^jy^kz^\ell(x^2- 1)(y^2-1)} 
\end{array}~:~ j+k+\ell = i-4 \right\} \\[8mm]
\begin{array}{c}
I^\otimes_i\Lambda^1(\square_3) \\
\hfill i\geq 1
\end{array} 
& 
\left\{~\qquad\qquad\qquad\qquad\qquad\qquad\qquad\textnormal{---\quad\textquotesingle\textquotesingle\quad---}\qquad\qquad\qquad\qquad\qquad\qquad\qquad~  ~:~ \max(k, j-1) = i - 1\right\}  \\[8mm]
\begin{array}{c}
\tilde I_i\Lambda^1(\square_3)\\
\hfill i\geq 4
\end{array}
 & \left\{\begin{array}{lll}
 \nfrdtab {y^{i-4}(y^2- 1)(z^2-1)}00 \\ 
 \nfrdtab {z^{i-4}(y^2- 1)(z^2-1)}00 \\ 
 \nfrdtab {0}{x^{i-4}(x^2- 1)(z^2-1)}0 \\
 \nfrdtab {0}{z^{i-4}(x^2- 1)(z^2-1)}0 \\
 \nfrdtab {0}{0}{x^{i-4}(x^2-1)(y^2-1)}\\
  \nfrdtab {0}{0}{y^{i-4}(x^2-1)(y^2-1)}\\
 \nfrdtab {x^jy^{i-j-4}(y^2- 1)(z^2-1)}{-x^{j-1}y^{i-j-3}(x^2- 1)(z^2-1)}0 \\
 \nfrdtab {x^jz^{i-j-4}(y^2- 1)(z^2-1)}{0}{-x^{j-1}z^{i-j-3}(x^2-1)(y^2-1)}\\
 \nfrdtab {0}{y^jz^{i-j-4}(x^2- 1)(z^2-1)}{-y^{j-1}z^{i-j-3}(x^2-1)(y^2-1)}~~\ast \\
 \hline
\nfrdtab {}{~:~1\leq j\leq i-4}{}\\
\multicolumn{3}{c}{\text{$\ast$: The last row is omitted when $i=5$ but included in all other cases. }}
 \end{array}\right\} \\[32mm]
\end{array}
$
\end{adjustwidth}

Note the single exception to the 1-form definitions that occurs in the case $\tilde I_5\Lambda^1(\square_3)$.
Without omitting the last row as indicated, the last three rows in the definition are linearly dependent.
For $i=4$, the last three rows are not included due to the vacuous condition $1\leq j\leq 0$ and for $i>5$, there is not enough repetition in the polynomial coefficients to cause a dependency.

\newpage

The 2-form spaces in 3D are constructed from the following spaces.
Each function in the sets listed below should be interpreted as the polynomial coefficient in the first column in front of $dydz$ plus the polynomial in the second column in front of $dxdz$ plus the polynomial in the third column in front of $dxdy$.\\
\begin{adjustwidth}{-.2in}{-.2in}
$
\begin{array}{rll}
\text{space} & \text{functions listed by $\left\{\begin{array}{lll} dydz & dxdz & dxdy \end{array}\right\}$} & \text{$i$ values} \\[2mm]
\hline\\[-2mm]
F^{\otimes}_i\Lambda^2(\square_3) & \left\{\begin{array}{lll}
 \nfrdtab {x^jy^kz^\ell(x \pm 1)}00 \\  
 \nfrdtab 0{x^ky^jz^\ell(y \pm 1)}0 \\ 
 \nfrdtab 00{x^ky^\ell z^j(z  \pm 1)} 
\end{array}~:~ \max(j+1, k, \ell) = i-1 \right\} & i\geq 1 \\[8mm]
F_i\Lambda^2(\square_3) & \left\{\begin{array}{lll}
 \nfrdtab {y^jz^k(x \pm 1)}00 \\  
 \nfrdtab 0{x^jz^k(y \pm 1)}0 \\ 
 \nfrdtab 00{x^jy^k(z  \pm 1)} 
\end{array}~:~ j+k = i \right\} & i\geq 0 \\[8mm]
\tilde F_i\Lambda^2(\square_3) & \left\{\begin{array}{lll}
 \nfrdtab {(i+1)y^i(x \pm 1)}{y^{i-1}(y^2-1)}0 \\  
 \nfrdtab {(i+1)z^i(x \pm 1)}0{-z^{i-1}(z^2-1)} \\  
 \nfrdtab {x^{i-1}(x^2-1)}{(i+1)x^i(y \pm 1)}0 \\ 
 \nfrdtab 0{(i+1)z^i(y \pm 1)}{z^{i-1}(z^2-1)} \\ 
 \nfrdtab {-x^{i-1}(x^2-1)}0{(i+1)x^i(z  \pm 1)} \\
 \nfrdtab 0{y^{i-1}(y^2-1)}{(i+1)y^i(z  \pm 1)} \\
 \nfrdtab {(i+2)y^{j}z^{i-j}(x\pm 1)}{y^{j-1}z^{i-j}(y^2-1)}{-y^jz^{i-j-1}(z^2- 1)} \\
 \nfrdtab {x^{j-1}z^{i-j}(x^2-1)}{(i+2)x^{j}z^{i-j}(y\pm 1)}{x^jz^{i-j-1}(z^2- 1)} \\
 \nfrdtab {-x^{j-1}y^{i-j}(x^2-1)}{x^{j}y^{i-j-1}(y^2-1)}{(i+2)x^jy^{i-j}(z \pm 1)} \\
 \hline
\nfrdtab {}{~:~1\leq j\leq i-1}{}\\
\end{array}
\right\} & i\geq 1 \\[20mm]
I_i\Lambda^2(\square_3) & \left\{\begin{array}{lll}
 \nfrdtab {x^jy^kz^\ell(x^2- 1)}00 \\  
 \nfrdtab 0{x^jy^kz^\ell(y^2-1)}0 \\ 
 \nfrdtab 00{x^jy^kz^\ell(z^2- 1)} 
\end{array}~:~ j+k+\ell = i-2 \right\} & i\geq 2 \\[8mm]
I^{\otimes}_i\Lambda^2(\square_3) & 
\left\{ ~\qquad\qquad\qquad\qquad\textnormal{---\quad\textquotesingle\textquotesingle\quad---}\qquad\qquad\qquad\qquad~ ~:~ \max(j, k-1, \ell-1) = i-1 \right\} & i\geq 1\\[8mm]
\tilde{I}_i\Lambda^2(\square_3) & \left\{\begin{array}{llll}
 \nfrdtab {x^{i-2}(x^2- 1)}00 \\  
 \nfrdtab 0{y^{i-2}(y^2-1)}0 \\ 
 \nfrdtab 00{z^{i-2}(z^2- 1)} \\
 \nfrdtab {x^{i-j-2}y^j(x^2- 1)}{-x^{i-j-1}y^{j-1}(y^2-1)}{0} & \ast \\
 \nfrdtab {x^{i-j-2}z^j(x^2- 1)}{0}{x^{i-j-1}z^{j-1}(z^2-1)} & \ast\\
 \nfrdtab {0}{y^{i-j-2}z^j(y^2- 1)}{-y^{i-j-1}z^{j-1}(z^2-1)} & \ast\\
 \nfrdtab {x^jy^kz^\ell(x^2-1)}{-x^{j+1}y^{k-1}z^\ell(y^2-1)}{x^{j+1}y^kz^{\ell-1}(z^2-1)} & \star\\
  \hline
  \multicolumn{3}{c}{
  \begin{array}{rl} 
  \text{for rows with $\ast$}: & ~1\leq j\leq i-2,\\
  \text{for row with $\star$}: & ~ k\geq 1, ~~\ell\geq 1, ~~j+k+\ell= i-2.
  \end{array}
  }\\
\end{array}
\right\} & i\geq 2 \\[8mm]
\end{array}
$
\end{adjustwidth}

\newpage
The 3-form spaces in 3D are constructed from the following spaces.
\[
\begin{array}{rll}
\text{space} & \text{functions} & \text{$i$ values} \\[2mm]
\hline\\[-2mm]
I_i\Lambda^3(\square_3) & \left\{x^jy^kz^\ell~dxdydz ~:~  j + k + \ell = i \right\} & i\geq 0 \\[8mm]
I^{\otimes}_i\Lambda^3(\square_3) & \left\{~\quad\textnormal{---\quad\textquotesingle\textquotesingle\quad---}\quad~ ~:~ max(j, k, \ell) = i-1\right\} & i\geq 1 \\[8mm]
\end{array}
\]

\subsection*{Subspace structure of computational bases.}
We have defined computational bases in terms of shared subspace structures.
By comparing these structures, we can better understand the nature of the various approximation spaces and devise efficient implementation strategies.

We first examine the 0-form bases in the $n=3$ case (the $n=2$ case is similar).
Comparing $\calQ_r^-\Lambda^0$ to $\calS_r\Lambda^0$, we see that the bases are identical for $m=0$ and 1 but differ for $m=2$ and 3.
For $m=2$, we examine the definitions for $F_i\Lambda^0(\square_3)$ and $F_i^\otimes\Lambda^0(\square_3)$ and see that the same kinds of functions are involved -- the difference is simply a matter of the constraints imposed on $j$ and $k$ in terms of $i$.
Thinking of a choice for $j$ and $k$ as an order pair in the plane, the constraint $j+k = i-4$ in the definition of  $F_i\Lambda^0(\square_3)$ allows pairs lying on a diagonal line in the first quadrant while the constraint $\max(j, k) = i-1$ in the definition of $F_i^\otimes\Lambda^0(\square_3)$ allows pairs lying on the top and right edges of a rectangle with its lower left corner at the origin.
Accordingly, when these spaces are summed to form computational bases, the result is a ``total degree'' type constraint for $\calS_r\Lambda^0$ and an ``each variable degree'' type constraint for $\calQ_r^-\Lambda^0$.
This type of growth is not a surprise -- the serendipity spaces are meant to provide approximation of a maximum total polynomial order while tensor product spaces are meant to provide approximation of a maximum order in each variable.
Note that since $\calS_r^-\Lambda^0=\calS_r\Lambda^0$, the computational bases are accordingly identical. 

Moving to the 1-form bases, we start with the $n=2$ case.
We note that $\bigoplus_{i=0}^{r-1}E_i\Lambda^1(\square_2)$ is a subspace of $\calQ_r^-\Lambda^1(\square_2)$, $\calS_r\Lambda^1(\square_2)$, and $\calS_r^-\Lambda^1(\square_2)$.
In the $\calQ^-\Lambda^1$ and $\calS^-\Lambda^1$ spaces, this is the complete portion of the basis for $m=1$, however, the serendipity subspace for $m=1$ has the additional direct summand $\tilde E_r\Lambda^1(\square_2)$.
Note that $\tilde E_r\Lambda^1(\square_2)$ consists of exactly four functions, one for edge of $\square_2$, and these functions take the form of something ``expected'' plus something that vanishes on the boundary.
For instance:
\[\underbrace{(r+1)x^r(y + 1)dx}_{\text{associated to edge $\{y=1\}$}} + \underbrace{x^{r-1}(x^2 - 1)dy}_{\text{zero trace on $\p\square_2$}} \in \tilde E_r\Lambda^1(\square_2).\]

The $m=2$ subspaces for $k=1$ are different for each of the three families.
For the $\calQ_r^-\Lambda^1$ space, the $m=2$ subspace is made of $F_i^\otimes\Lambda^1$ spaces, which are the kinds of bases typically used to implement Raviart-Thomas spaces.
For the serendipity spaces, the $m=2$ subspace is comprised only of $F_i\Lambda^1$ spaces, which differ from $F_i^\otimes\Lambda^1$ only in the constraint on $j$ and $k$ in terms of $i$ (similar to the $m=1$ case).
The $m=2$ subspace of the trimmed serendipity spaces is slightly different in that it replaces the top order space $F_r\Lambda^1(\square_2)$ by $\tilde F_r\Lambda^1(\square_2)$.
The first two functions listed in the definition of $\tilde F_r\Lambda^1(\square_2)$ fit with the pattern established by the third if we allow $k=0$ and $k=i-1$ and interpret $x^{-1}$ and $y^{-1}$ as 0.  
The difference between the trimmed and non-trimmed spaces can be explained by the containment properties $\calP_r\Lambda^1\subset\calS_r\Lambda^1$ while $\calP_{r-1}\Lambda^1\subset\calP_r^-\Lambda^1\subset\calS_r^-\Lambda^1$ (see~\cite{GK2016}).
We thus expect that the ``top degree'' functions in the trimmed serendipity space will differ slightly from those of the serendipity spaces, as is indeed the case with these bases.

Even more subtle distinctions among the three 1-form spaces become apparent in the $n=3$ case.
The $m=1$ subspaces are direct analogues of the $n=2$ case, with a more involved definition of the $\tilde E_i\Lambda^1$ space due to the presence of a third variable.
Here, there are exactly 12 functions, each associated to a distinct edge of $[-1,1]^3$.
A detailed discussion of the properties of the functions in $\tilde E_1\Lambda^1(\square_3)$ was given as part of the introduction in Section~\ref{sec:intro}; higher orders of $r$ can be explained in the same way.
Note that for any $f\prec\square_3$ with $\dim f=2$, we have
\begin{align}
\tr_{f} E_i\Lambda^1(\square_3) & =E_i\Lambda^1(f)\\
\tr_{f}\tilde E_i\Lambda^1(\square_3) & =\tilde E_i\Lambda^1(f)
\end{align}
meaning the computational bases respect the geometry for the edge spaces.

Similar trace properties hold for the $m=2$ subspaces of the 1-form spaces. 
For any $f\prec\square_3$ with $\dim f=2$, we have
\begin{align}
\tr_{f} F_i\Lambda^1(\square_3) &=F_i\Lambda^1(f) \\
\tr_{f} \hat F_r\Lambda^1(\square_3) &=F_r\Lambda^1(f) \\
\tr_{f} \tilde F_r\Lambda^1(\square_3) &= \tilde F_r\Lambda^1(f)
\end{align}
These trace properties help us understand the difference between the $\hat F$ and $\tilde F$ spaces. 
Excluding the six functions that the two spaces share, each $\hat F$ function is ``something that vanishes on all $\p\square_3$ except a particular face'' plus ``something that vanishes on all $\p\square_3$,'' whereas each $\tilde F$ function is ``something that looks like an element of $\tilde F$ on a particular face''.
For instance, taking $i=3$, $j=1$, we compare
\begin{align}
\underbrace{4x(z + 1)(y^2-1) dx}_{\text{associated to face $\{z=1\}$}} &  + \underbrace{(x^2-1)(y^2-1)dz}_{\text{zero trace on $\p\square_3$}} & \in \hat F_3\Lambda^1(\square_3) \\
\underbrace{x(z + 1)(y^2-1) dx}_{\text{associated to face $\{z=1\}$}} &  +\underbrace{(-1) y(z+ 1)(x^2-1) dy}_{\text{required for $\tilde F_3\Lambda^1(\square_2)$ on $\{z=1\}$}} & \in \tilde F_3\Lambda^1(\square_3)
\end{align}
This comparison highlights the important and surprising fact that functions in $\hat F_3\Lambda^1(\square_3)$ contribute to the approximation power on the interior of the cube through terms such as $(x^2-1)(y^2-1)dz$.

Next we turn to the 2-form bases in the $n=3$ case, focusing on $\calS_r\Lambda^2(\square_3)$ and $\calS^-_r\Lambda^2(\square_3)$.
First, consider the overall subspace decompositions:
\begin{align}
\calS_r\Lambda^2(\square_3) & =  \left(\ds\bigoplus_{i=0}^{r-1} F_i\Lambda^2(\square_3) \oplus \tilde F_r\Lambda^2(\square_3)\right)\oplus \ds\bigoplus_{i=2}^r I_i\Lambda^3(\square_3)  \\
\calS^-_r\Lambda^2(\square_3) & =  \ds\bigoplus_{i=0}^{r-1} F_i\Lambda^2(\square_3) \oplus \left(\ds\bigoplus_{i=2}^{r-1} I_i\Lambda^3(\square_3) \oplus \tilde I_r\Lambda^2(\square_3) \right)
\end{align}
Observe that $\calS_r\Lambda^2(\square_3)$ has a special definition of functions associated to the top order on faces of co-dimension 1 while $\calS^-\Lambda^2(\square_3)$ has a special definition of functions associated to the top order on the interior of $\square_3$.
An analogous statement could be made about $\calS_r\Lambda^1(\square_2)$ and $\calS^-_r\Lambda^1(\square_2)$ as they have the same subspace structure with everything shifted down by a dimension:
\begin{align}
\calS_r\Lambda^1(\square_2) & =  \left(\ds\bigoplus_{i=0}^{r-1} E_i\Lambda^1(\square_2) \oplus \tilde E_r\Lambda^1(\square_2)\right)\oplus \ds\bigoplus_{i=2}^r F_i\Lambda^2(\square_2)  \\
\calS^-_r\Lambda^1(\square_2) & =  \ds\bigoplus_{i=0}^{r-1} E_i\Lambda^1(\square_2) \oplus \left(\ds\bigoplus_{i=2}^{r-1} F_i\Lambda^2(\square_2) \oplus \tilde F_r\Lambda^1(\square_2) \right)
\end{align}
Further, each function in $\tilde F_i\Lambda^2(\square_3)$ is associated to a face of $[-1,1]^3$ according to where a factor of the form $(\ast \pm 1)$ appears.
For instance,
\begin{equation}
\underbrace{4yz(x + 1)~dydz}_{\text{associated to face $\{x=1\}$}} +
\underbrace{z(y^2-1)~dxdx}_{\text{zero trace on $\p\square_3$}} -
\underbrace{y(z^2- 1)~dxdy}_{\text{zero trace on $\p\square_3$}} ~\in \tilde F_2\Lambda^2(\square_3).
\end{equation}

Looking at the $I_i\Lambda^2(\square_3)$ and $\tilde I_i\Lambda^2(\square_3)$ spaces, note that all functions in these spaces must have zero trace on $\p\square_3$, forcing a factor of $(x^2-1)$ in front of $dydz$ terms and so forth with the other alternators. 
The definition explains how monomials of total degree $i-2$ may appear in the polynomial coefficients.
As in the $\tilde F_r\Lambda^1(\square_2)$ subspaces, the definition could be simplified if we interpret terms like $x^{-1}yz$ as 0, but we have again opted for clarity of definition over conciseness.

Finally, we note that the top order form spaces, namely $k=2$ for $n=2$ and $k=3$ for $n=3$, are exactly what is expected: order $r$ tensor product spaces for $\calQ^-$, order $r$ polynomial spaces for $\calS$, and order $r-1$ polynomial spaces for $\calS^-$.

\subsection*{Implementation of computational bases.}
The subspace structure allows efficient implementation of the computational bases. 
For instance, we can construct the rather involved 1-form spaces for both serendipity and trimmed serendipity elements of any order $r$ just by coding seven subspace definitions: $E_i\Lambda^1$, $\tilde E_i\Lambda^1$,  $F_i\Lambda^1$,  $\hat F_i\Lambda^1$,  $\tilde F_i\Lambda^1$,  $I_i\Lambda^1$, and $\tilde I_i\Lambda^1$, followed by an assembly command.
In the next section, we explain how we did exactly this using \textit{SageMath}.

\section{SageMath code for basis construction and verification}
\label{sec:bases-sage}

We have written code in the open source software system \textit{SageMath}~\cite{sagemath} to generate lists of computational bases and to verify (by a simple linear algebra trick) that the lists are in fact bases of the applicable polynomial differential form space.
We chose to use \textit{SageMath} for this purpose as it has a package specifically for storing and manipulating polynomial differential forms while also allowing import of \texttt{numpy} and other open source packages.
The code for $n=3$ is available at \url{http://math.arizona.edu/~agillette/research/computationalBases/}, with some highlights described below.

\subsection*{Code for the Koszul operator.}
The exterior derivative operation ($d$, from~\eqref{eq:def-d}) is a built-in function of the \texttt{DifferentialForms} class of \textit{SageMath}, but the Koszul operator, $\kappa$, is not.
We wrote a routine for $\kappa$, based on \eqref{eq:def-kappa}:\\[-2mm]

\begin{verbatim}
def k(form):
    form._cleanup()
    kappa_form = DifferentialForm(F, form._degree - 1)
    for comp in form._components:
        fun = form._components[comp].factor()
        i = 1
        for n, coord in enumerate(F.base_space().coordinates()):
            if n in comp:
                newcomp = tuple(a for a in comp if a != n)
                kappa_form[newcomp] += fun*coord*(-1)^(i+1)
                i += 1
    kappa_form._cleanup()
    return kappa_form
\end{verbatim}
Given an input \texttt{form}, the function \texttt{k} constructs an output \texttt{kappa$\mathunderscore$form} of one form order less by iterating through each alternator appearing in \texttt{form}.
For each alternator in the input (\texttt{comp}), it checks each coordinate (\texttt{n}) from the globally fixed base space of variables (e.g.\ `$x,y,z$'), to see if if the coordinate is present in the alternator.
If so, a new component (\texttt{newcomp}) is added to the output, based on the formula~\eqref{eq:def-kappa}.
The \texttt{cleanup} command removes empty components before the output is returned.

\subsection*{Code for the $\calH_r\Lambda^k(\bbR^n)$ spaces.}
Recall that a basis for $\calH_r\Lambda^k(\bbR^n)$, the space of homogeneous $k$-forms of polynomial degree $r$, is the set of form monomials $x^\alpha dx_\sigma$ such that $|\alpha|=r$ and $|\sigma|=k$.
This basis is built as follows:

\begin{verbatim}
def HrLkRn(r,k,n):
    U = CoordinatePatch(var(n))
    F = DifferentialForms(U)
    monomials = []
    for mon in itertools.product(var(n),repeat=r):
        monomials.append(prod(mon))
    allforms = []
    for dxs in itertools.combinations(range(0,F.ngens()),k):
        for mon in monomials:
            form = DifferentialForm(F,k)
            form[dxs] = mon
            if form not in allforms:
                allforms.append(form)
    return allforms
\end{verbatim}

The inputs $r$ and $k$ should be non-negative integers but $n$ should be a string variables to be used, separated by commas, e.g. \texttt{'x,y,z'}.
The monomials of degree exactly $r$ are built with \texttt{itertools.product} and a list of all $k$-combinations of $n$ variables is collected with \texttt{itertools.combinations}.
The result of the nested \texttt{for} loops is a list of length (\# of monomials)$\cdot$(\# of alternators), containing exactly the desired basis for $\calH_r\Lambda^k(\bbR^n)$.

\subsection*{Code for the $\calJ_r\Lambda^k(\bbR^n)$ spaces.}
We use the \texttt{k} and \texttt{HrLkRn} routines to build a basis for $\calJ_r \Lambda^k(\bbR^n)$.  
The code is based on the definition in~\eqref{eq:jrk-def}, but with the $\kappa$ operator pulled out:
\[\calJ_r \Lambda^k(\bbR^n) = \kappa \sum_{l \geq 1} \calH_{r+l-1, l} \Lambda^{k+1}(\bbR^n)\]

\begin{verbatim}
def JrLkRn(r,j,n):
    Hs = []
    for l in range(1,len(var(n))-j):
        eH = HrLkRn(r+l-1,j+1,n)
        if j+1 != 0:
            for form in eH:
                for comp in form._components:
                    ldeg = 0
                    fun = form[comp]
                    for w,coord in enumerate(F.base_space().coordinates()):
                        if w not in comp and fun.degree(coord) == 1:    
                            ldeg += 1
                    if ldeg >= l:
                        Hs.append(form)
    final = []
    for form1 in Hs:
        new = k(form1)
        if new not in final:
            final.append(new)
    return final
\end{verbatim}

The inputs are again non-negative integers $r$ and $j$ (a local variable replacing $k$, to avoid confusion with the $\texttt{k}$ command), and a string variables for $n$, as in the routine for \texttt{HrLkRn}.
For each possible $\ell$ value (typewriter font version is \texttt{l}), starting at $\ell=1$, the code builds the basis for $\calH_{r+\ell-1}\Lambda^{j+1}(\bbR^n)$.
Each basis element is checked to see if it has linear degree at least $\ell$. 
If so, it is added to a temporary list \texttt{Hs}. 
When the for loop is done, each element in \texttt{Hs} has the $\kappa$ operator applied to it and is appended (avoiding duplicates but not checking for linear dependencies) to the final list that is returned.

\subsection*{Assembling standard spanning sets and computational bases.}
Using the routines for \texttt{HrLkRn}, \texttt{JrLkRn}, and similar ideas, we create `standard' spanning sets for the five major polynomial differential form spaces, based on their definitions in~\eqref{eq:prk-def} -- \eqref{eq:srmk-def}.
Note that while a basis is typically produced, some low-value $r$ and $k$ choices yield a linearly dependent set.
For instance, the code generates $d\calJ_1\Lambda^1(\bbR^3)$ as the span of
\[\left\{\begin{array}{l} -x~dy dz + y~dx dz + 2z~dx dy, \\ ~x~dy dz + 2y~dx dz + z~dx dy, \\ 2x~dy dz  + y~dx dz -z~dx dy\end{array}\right\},\]
but $\dim d\calJ_1\Lambda^1(\bbR^3)=2$.
Rather than arbitrarily remove one of these elements, we leave them all as a spanning set as it will not affect the subsequent constructions.

We use separate routines to construct the subspaces $V\Lambda^k$, $E_i\Lambda^k$, $F_i\Lambda^k$, $I_i\Lambda^k$ and their variants, as defined in Section~\ref{sec:bases-theory}.
These pieces are used to build computational bases, based on the definitions in Tables \ref{tab:dim2} and \ref{tab:dim3}.
For instance, a basis for $\calS_3\Lambda^1(\square_2)$ defined by the decomposition~(\ref{eq:S312-comp-bas-ex}) is assembled by the command \texttt{Spr13} (the \texttt{p} is for `plus' as opposed to `minus' in the trimmed case), as in this \textit{SageMath} session:
\begin{verbatim}
  sage: attach("construct-tools-n3.sage")
  sage: Spr13(2)
  [(y + 1)*(z + 1)*dx,
   ...
   (y^2 - 1)*(x - 1)*dz]
\end{verbatim}

\subsection*{Basis verification procedure.}
In Section~\ref{sec:bases-theory}, we explained the geometric association of the elements in the various bases, but it remains to verify that the stated definitions are indeed bases for the relevant spaces.
We devise a computational procedure based on the following simple linear algebra results.

\begin{lemma}
Let $V$ be a finite dimensional, real vector space. Let $A$ and $B$ be sets of vectors in $V$. If 
\[ \dim~\myspan_V(A) = \dim~\myspan_V(A \cup B) = \dim~\myspan_V(B), \]
then $\myspan_V(A) = \myspan_V(B)$.
\end{lemma}

\begin{proof}
Since $V$ is finite dimensional, there exists $N=\dim~\myspan_V(A\cup B)$, $0\leq N<\infty$.
Since $A\subset A\cup B$, $\dim~\myspan_V(A) \leq N$, but by hypothesis, this is an equality.
The only $N$-dimensional subspace of an $N$-dimensional space is the space itself, whereby $\myspan_V(A) = \myspan_V(A\cup B)$.
By the same argument, $\myspan_V(B) = \myspan_V(A\cup B)$ and the conclusion follows.
\end{proof}

\begin{corollary}
Let $A$ be a set of polynomial differential $k$-forms such that $\myspan(A)=\mathcal X\Lambda^k(\square_n)$.
Let $B$ be a set of polynomial differential $k$-forms such that 
\[\#B = \dim~\myspan(B)=\dim~\myspan (A \cup B) =\dim~\myspan(A).\]
Then $B$ is a basis for $\mathcal X\Lambda^k(\square_n)$.
\end{corollary}

\begin{proof}
The first equality ensures that $B$ is a set of linearly independent elements.
The second and third equalities ensure the hypothesis of the Lemma is satisfied.
It follows that $B$ is a spanning set for $\mathcal X\Lambda^k(\square_n)$ and hence a basis for it.
\end{proof}

In light of the Corollary, our procedure for checking that our stated lists are bases is as follows:\\

\noindent\textit{Basis Verification Algorithm}
\begin{enumerate}[1.]
\item Create $A$, the standard spanning set for a polynomial differential form space.\\[-7mm]
\item Create $B$, the claimed computational basis set for the same space.\\[-7mm]
\item Expand each element of each set in terms of form monomials.  Store coefficients of these expansions in matrices $\mathbb A$ and $\mathbb B$, respectively.\\[-7mm]
\item Define matrix $\mathbb C$ by taking $\mathbb A$ and adding the rows of $\mathbb B$ below it.\\[-7mm]
\item Compute the ranks of $\mathbb A$, $\mathbb B$ and $\mathbb C$.  If they all agree, then $B$ is a basis for the space spanned by $A$.
\end{enumerate}

The code included with the manuscript includes a command called \texttt{basis\textunderscore  check()}, which carries out the Basis Verification Algorithm for $\calS_1\Lambda^1(\square_3)$.
The output of this command includes the computational basis for $\calS_1\Lambda^1(\square_3)$, a statement of its length (24), and a verification that the three matrices constructed all have rank 24, as desired.
The interested reader can recreate this output by the following commands:
\begin{verbatim}
  sage: attach("construct-tools-n3.sage")
  sage: basis_check()
\end{verbatim}
On a 2.9 GHz Mac Desktop, this command took about 1 minute to run.
We tested many more cases on a department-owned Xeon processor with 32 GB of RAM.

Step 3 of the algorithm is the only part of the coding that introduces any subtleties.
Here is a piece of the code that builds the list of coefficients for a set of $1$-forms when $n=3$:
\begin{verbatim}
MR3.<x,y,z>=QQbar[];MR3
...
def find_coeffs_k1(basis_list,deg):
    final_list=[]
    for form in basis_list:
        coeff_list=[]   
        for i in range(0,3): 
            for xdeg in range(0,deg+1):
                for ydeg in range(0,deg+1):
                    for zdeg in range(0,deg+1):
                        poly=MR3(form[i])
                        c=poly.coefficient({x:xdeg,y:ydeg,z:zdeg})
                        coeff_list.append(c)
        final_list.append(coeff_list)
    return final_list
\end{verbatim}
The inputs are a list of differential 1-forms and a degree, which should be at least as big as the maximum exponent appearing in any polynomial coefficient of the input list.
The command \texttt{MR3(form[i])} tells \textit{SageMath} to interpret the $i$th component of \texttt{form} as an element of a multivariate polynomial ring in $x$, $y$, $z$ over an algebraic field. 
This `typing' as a polynomial allows access to the object property \texttt{coefficient} that is called in the next line.
For instance, the form $(5x+4yz+3y+2z+1) dx$ with $\texttt{deg}=1$ yields the list of coefficients $[\begin{array}{cccccccc} 1 & 2 & 3 & 4 & 5 & 0 & \cdots & 0\end{array}]$, a list of 24 numbers since there are 24 form monomials $x^{\alpha}dx_{\sigma}$ with $\sigma\in\{1,2,3\}$ satisfying $0\leq \alpha_i\leq 1$.

A small additional issue arises in the $k=0$ case. 
The \texttt{DifferentialForms} package treats coefficients of $k$-forms for $k>0$ as elements of a ``Symbolic Ring'' class, but it treats 0-forms as elements of an ``Algebra of differential forms in the variables x, y, z.''
Thus to re-type our 0-form lists as multivariate polynomials, we have \textit{Sagemath} re-read the forms as strings, which reclassifies them as elements of a ``Symbolic Ring.''
For example, in the definition of \texttt{find\textunderscore coeffs\textunderscore k0}, we have:
\begin{verbatim}
basis_list=sage_eval(str(basis_list),locals={'x':x, 'y':y, 'z':z})
\end{verbatim}
Beyond these minor issues, we found the \texttt{DifferentialForms} package to be very well implemented and user friendly.

\section{Conclusion}
\label{sec:conc}
We have presented a complete list of computational basis functions for $\calQ_r^-\Lambda^k(\square_n)$, $\calS_r\Lambda^k(\square_n)$, and $\calS_r^-\Lambda^k(\square_n)$, for $n=2$ and 3, $0\leq k\leq n$, and $r\geq 1$.
As is evident from the subspace structure of these bases seen in Tables~\ref{tab:dim2} and \ref{tab:dim3}, all the tensor product bases and the $k=0$ and $k=n$ serendipity bases are naturally \textit{hierarchical}, meaning the basis of order $r+1$ can be constructed by simply adding elements to the basis of of order $r$.
The serendipity bases for $k=1$ and $k=2$ (for $n=3$) are not hierarchical due to the presence of subspaces like $\tilde E_r$ occurring in the top order polynomial degree.
This special treatment of top order is a feature of serendipity spaces, not a bug, and thus is not avoidable in basis construction.
We suspect that the bases defined here are ``as hierarchical as possible.''
By presenting these bases in terms of differential forms, we hope to aid their implementation into the Unified Form Language~\cite{ALORW2014} used by FEniCS and other open source, multi-purpose finite element packages.

\bibliography{bib-koszul}{}
\bibliographystyle{abbrv}

\begin{appendices}

\section{Lists of $\calS_r\Lambda^k$ and $\calS_r^-\Lambda^k$ basis functions}
\label{app:lists}

Below are computational basis functions for the serendipity space $\calS_r\Lambda^k(\square_3)$ and the trimmed serendipity space $\calS_r^-\Lambda^k(\square_3)$ for $r=1$ to 3, $k=0$ to 2, and $n=3$.
To the best of our knowledge, for $k=1$ and $k=2$, these bases have not appeared in the literature previously.
All bases were generated using \textit{SageMath} code and verified using the Basis Verification Algorithm described in the paper.
The code is available at \url{http://math.arizona.edu/~agillette/research/computationalBases/}.

\begin{tiny}
\[
\begin{array}{cl}
\calS_1\Lambda^0(\square_3),~~\calS_1^-\Lambda^0(\square_3) 
&
\begin{array}{lll}
\hline
(x + 1)(y + 1)(z + 1) \\
(x + 1)(y + 1)(z - 1) \\
(x + 1)(y - 1)(z + 1) \\
(x + 1)(y - 1)(z - 1) \\
(x - 1)(y + 1)(z + 1) \\
(x - 1)(y + 1)(z - 1) \\
(x - 1)(y - 1)(z + 1) \\
(x - 1)(y - 1)(z - 1) \\
\end{array} 
\end{array} \]
\[
\begin{array}{cl}
\calS_2\Lambda^0(\square_3),~~\calS_2^-\Lambda^0(\square_3)
&
\begin{array}{lll}
\hline
(x + 1)(y + 1)(z + 1) \\
(x + 1)(y + 1)(z - 1) \\
(x + 1)(y - 1)(z + 1) \\
(x + 1)(y - 1)(z - 1) \\
(x - 1)(y + 1)(z + 1) \\
(x - 1)(y + 1)(z - 1) \\
(x - 1)(y - 1)(z + 1) \\
(x - 1)(y - 1)(z - 1) \\
(x^2 - 1)(y + 1)(z + 1) \\
(x^2 - 1)(y + 1)(z - 1) \\
(x^2 - 1)(y - 1)(z + 1) \\
(x^2 - 1)(y - 1)(z - 1) \\
(y^2 - 1)(x + 1)(z + 1) \\
(y^2 - 1)(x - 1)(z + 1) \\
(y^2 - 1)(x + 1)(z - 1) \\
(y^2 - 1)(x - 1)(z - 1) \\
(z^2 - 1)(x + 1)(y + 1) \\
(z^2 - 1)(x + 1)(y - 1) \\
(z^2 - 1)(x - 1)(y + 1) \\
(z^2 - 1)(x - 1)(y - 1) \\

\end{array} 
\end{array} \]
\[
\begin{array}{cl}
\calS_3\Lambda^0(\square_3),~~\calS_3^-\Lambda^0(\square_3)
&
\begin{array}{lll}
\hline
(x + 1)(y + 1)(z + 1) \\
(x + 1)(y + 1)(z - 1) \\
(x + 1)(y - 1)(z + 1) \\
(x + 1)(y - 1)(z - 1) \\
(x - 1)(y + 1)(z + 1) \\
(x - 1)(y + 1)(z - 1) \\
(x - 1)(y - 1)(z + 1) \\
(x - 1)(y - 1)(z - 1) \\
(x^2 - 1)(y + 1)(z + 1) \\
(x^2 - 1)(y + 1)(z - 1) \\
(x^2 - 1)(y - 1)(z + 1) \\
(x^2 - 1)(y - 1)(z - 1) \\
(y^2 - 1)(x + 1)(z + 1) \\
(y^2 - 1)(x - 1)(z + 1) \\
(y^2 - 1)(x + 1)(z - 1) \\
(y^2 - 1)(x - 1)(z - 1) \\
(z^2 - 1)(x + 1)(y + 1) \\
(z^2 - 1)(x + 1)(y - 1) \\
(z^2 - 1)(x - 1)(y + 1) \\
(z^2 - 1)(x - 1)(y - 1) \\
(x^2 - 1)x(y + 1)(z + 1) \\
(x^2 - 1)x(y + 1)(z - 1) \\
(x^2 - 1)x(y - 1)(z + 1) \\
(x^2 - 1)x(y - 1)(z - 1) \\
(y^2 - 1)(x + 1)y(z + 1) \\
(y^2 - 1)(x - 1)y(z + 1) \\
(y^2 - 1)(x + 1)y(z - 1) \\
(y^2 - 1)(x - 1)y(z - 1) \\
(z^2 - 1)(x + 1)(y + 1)z \\
(z^2 - 1)(x + 1)(y - 1)z \\
(z^2 - 1)(x - 1)(y + 1)z \\
(z^2 - 1)(x - 1)(y - 1)z \\

\end{array} 
\end{array} \]
\[
\begin{array}{cl}
\calS_1\Lambda^1(\square_3)
&
\begin{array}{lll}
dx & dy & dz \\
\hline
(y + 1)(z + 1) & 0 & 0 \\
(y + 1)(z - 1) & 0 & 0 \\
(y - 1)(z + 1) & 0 & 0 \\
(y - 1)(z - 1) & 0 & 0 \\
0 & (x + 1)(z + 1) & 0 \\
0 & (x - 1)(z + 1) & 0 \\
0 & (x + 1)(z - 1) & 0 \\
0 & (x - 1)(z - 1) & 0 \\
0 & 0 & (x + 1)(y + 1) \\
0 & 0 & (x + 1)(y - 1) \\
0 & 0 & (x - 1)(y + 1) \\
0 & 0 & (x - 1)(y - 1) \\
2x(y + 1)(z + 1) & (x^2 - 1)(z + 1) & (x^2 - 1)(y + 1) \\
2x(y + 1)(z - 1) & (x^2 - 1)(z - 1) & (x^2 - 1)(y + 1) \\
2x(y - 1)(z + 1) & (x^2 - 1)(z + 1) & (x^2 - 1)(y - 1) \\
2x(y - 1)(z - 1) & (x^2 - 1)(z - 1) & (x^2 - 1)(y - 1) \\
(y^2 - 1)(z + 1) & 2(x + 1)y(z + 1) & (y^2 - 1)(x + 1) \\
(y^2 - 1)(z + 1) & 2(x - 1)y(z + 1) & (y^2 - 1)(x - 1) \\
(y^2 - 1)(z - 1) & 2(x + 1)y(z - 1) & (y^2 - 1)(x + 1) \\
(y^2 - 1)(z - 1) & 2(x - 1)y(z - 1) & (y^2 - 1)(x - 1) \\
(z^2 - 1)(y + 1) & (z^2 - 1)(x + 1) & 2(x + 1)(y + 1)z \\
(z^2 - 1)(y - 1) & (z^2 - 1)(x + 1) & 2(x + 1)(y - 1)z \\
(z^2 - 1)(y + 1) & (z^2 - 1)(x - 1) & 2(x - 1)(y + 1)z \\
(z^2 - 1)(y - 1) & (z^2 - 1)(x - 1) & 2(x - 1)(y - 1)z \\
\end{array} 
\end{array} \]
\[
\begin{array}{cl}
\calS_2\Lambda^1(\square_3)
&
\begin{array}{lll}
dx & dy & dz \\
\hline
(y + 1)(z + 1) & 0 & 0 \\
(y + 1)(z - 1) & 0 & 0 \\
(y - 1)(z + 1) & 0 & 0 \\
(y - 1)(z - 1) & 0 & 0 \\
0 & (x + 1)(z + 1) & 0 \\
0 & (x - 1)(z + 1) & 0 \\
0 & (x + 1)(z - 1) & 0 \\
0 & (x - 1)(z - 1) & 0 \\
0 & 0 & (x + 1)(y + 1) \\
0 & 0 & (x + 1)(y - 1) \\
0 & 0 & (x - 1)(y + 1) \\
0 & 0 & (x - 1)(y - 1) \\
x(y + 1)(z + 1) & 0 & 0 \\
x(y + 1)(z - 1) & 0 & 0 \\
x(y - 1)(z + 1) & 0 & 0 \\
x(y - 1)(z - 1) & 0 & 0 \\
0 & (x + 1)y(z + 1) & 0 \\
0 & (x - 1)y(z + 1) & 0 \\
0 & (x + 1)y(z - 1) & 0 \\
0 & (x - 1)y(z - 1) & 0 \\
0 & 0 & (x + 1)(y + 1)z \\
0 & 0 & (x + 1)(y - 1)z \\
0 & 0 & (x - 1)(y + 1)z \\
0 & 0 & (x - 1)(y - 1)z \\
3x^2(y + 1)(z + 1) & (x^2 - 1)x(z + 1) & (x^2 - 1)x(y + 1) \\
3x^2(y + 1)(z - 1) & (x^2 - 1)x(z - 1) & (x^2 - 1)x(y + 1) \\
3x^2(y - 1)(z + 1) & (x^2 - 1)x(z + 1) & (x^2 - 1)x(y - 1) \\
3x^2(y - 1)(z - 1) & (x^2 - 1)x(z - 1) & (x^2 - 1)x(y - 1) \\
(y^2 - 1)y(z + 1) & 3(x + 1)y^2(z + 1) & (y^2 - 1)(x + 1)y \\
(y^2 - 1)y(z + 1) & 3(x - 1)y^2(z + 1) & (y^2 - 1)(x - 1)y \\
(y^2 - 1)y(z - 1) & 3(x + 1)y^2(z - 1) & (y^2 - 1)(x + 1)y \\
(y^2 - 1)y(z - 1) & 3(x - 1)y^2(z - 1) & (y^2 - 1)(x - 1)y \\
(z^2 - 1)(y + 1)z & (z^2 - 1)(x + 1)z & 3(x + 1)(y + 1)z^2 \\
(z^2 - 1)(y - 1)z & (z^2 - 1)(x + 1)z & 3(x + 1)(y - 1)z^2 \\
(z^2 - 1)(y + 1)z & (z^2 - 1)(x - 1)z & 3(x - 1)(y + 1)z^2 \\
(z^2 - 1)(y - 1)z & (z^2 - 1)(x - 1)z & 3(x - 1)(y - 1)z^2 \\
(y^2 - 1)(z + 1) & 0 & 0 \\
(y^2 - 1)(z - 1) & 0 & 0 \\
(z^2 - 1)(y + 1) & 0 & 0 \\
(z^2 - 1)(y - 1) & 0 & 0 \\
0 & (x^2 - 1)(z + 1) & 0 \\
0 & (x^2 - 1)(z - 1) & 0 \\
0 & (z^2 - 1)(x + 1) & 0 \\
0 & (z^2 - 1)(x - 1) & 0 \\
0 & 0 & (x^2 - 1)(y + 1) \\
0 & 0 & (x^2 - 1)(y - 1) \\
0 & 0 & (y^2 - 1)(x + 1) \\
0 & 0 & (y^2 - 1)(x - 1) \\
\end{array} 
\end{array} \]
\[
\begin{array}{cl}
\calS_3\Lambda^1(\square_3)
&
\begin{array}{lll}
dx & dy & dz \\
\hline
(y + 1)(z + 1) & 0 & 0 \\
(y + 1)(z - 1) & 0 & 0 \\
(y - 1)(z + 1) & 0 & 0 \\
(y - 1)(z - 1) & 0 & 0 \\
0 & (x + 1)(z + 1) & 0 \\
0 & (x - 1)(z + 1) & 0 \\
0 & (x + 1)(z - 1) & 0 \\
0 & (x - 1)(z - 1) & 0 \\
0 & 0 & (x + 1)(y + 1) \\
0 & 0 & (x + 1)(y - 1) \\
0 & 0 & (x - 1)(y + 1) \\
0 & 0 & (x - 1)(y - 1) \\
x(y + 1)(z + 1) & 0 & 0 \\
x(y + 1)(z - 1) & 0 & 0 \\
x(y - 1)(z + 1) & 0 & 0 \\
x(y - 1)(z - 1) & 0 & 0 \\
0 & (x + 1)y(z + 1) & 0 \\
0 & (x - 1)y(z + 1) & 0 \\
0 & (x + 1)y(z - 1) & 0 \\
0 & (x - 1)y(z - 1) & 0 \\
0 & 0 & (x + 1)(y + 1)z \\
0 & 0 & (x + 1)(y - 1)z \\
0 & 0 & (x - 1)(y + 1)z \\
0 & 0 & (x - 1)(y - 1)z \\
x^2(y + 1)(z + 1) & 0 & 0 \\
x^2(y + 1)(z - 1) & 0 & 0 \\
x^2(y - 1)(z + 1) & 0 & 0 \\
x^2(y - 1)(z - 1) & 0 & 0 \\
0 & (x + 1)y^2(z + 1) & 0 \\
0 & (x - 1)y^2(z + 1) & 0 \\
0 & (x + 1)y^2(z - 1) & 0 \\
0 & (x - 1)y^2(z - 1) & 0 \\
0 & 0 & (x + 1)(y + 1)z^2 \\
0 & 0 & (x + 1)(y - 1)z^2 \\
0 & 0 & (x - 1)(y + 1)z^2 \\
0 & 0 & (x - 1)(y - 1)z^2 \\
4x^3(y + 1)(z + 1) & (x^2 - 1)x^2(z + 1) & (x^2 - 1)x^2(y + 1) \\
4x^3(y + 1)(z - 1) & (x^2 - 1)x^2(z - 1) & (x^2 - 1)x^2(y + 1) \\
4x^3(y - 1)(z + 1) & (x^2 - 1)x^2(z + 1) & (x^2 - 1)x^2(y - 1) \\
4x^3(y - 1)(z - 1) & (x^2 - 1)x^2(z - 1) & (x^2 - 1)x^2(y - 1) \\
(y^2 - 1)y^2(z + 1) & 4(x + 1)y^3(z + 1) & (y^2 - 1)(x + 1)y^2 \\
(y^2 - 1)y^2(z + 1) & 4(x - 1)y^3(z + 1) & (y^2 - 1)(x - 1)y^2 \\
(y^2 - 1)y^2(z - 1) & 4(x + 1)y^3(z - 1) & (y^2 - 1)(x + 1)y^2 \\
(y^2 - 1)y^2(z - 1) & 4(x - 1)y^3(z - 1) & (y^2 - 1)(x - 1)y^2 \\
(z^2 - 1)(y + 1)z^2 & (z^2 - 1)(x + 1)z^2 & 4(x + 1)(y + 1)z^3 \\
(z^2 - 1)(y - 1)z^2 & (z^2 - 1)(x + 1)z^2 & 4(x + 1)(y - 1)z^3 \\
(z^2 - 1)(y + 1)z^2 & (z^2 - 1)(x - 1)z^2 & 4(x - 1)(y + 1)z^3 \\
(z^2 - 1)(y - 1)z^2 & (z^2 - 1)(x - 1)z^2 & 4(x - 1)(y - 1)z^3 \\
(y^2 - 1)(z + 1) & 0 & 0 \\
(y^2 - 1)(z - 1) & 0 & 0 \\
(z^2 - 1)(y + 1) & 0 & 0 \\
(z^2 - 1)(y - 1) & 0 & 0 \\
0 & (z^2 - 1)(x + 1) & 0 \\
0 & (z^2 - 1)(x - 1) & 0 \\
0 & (x^2 - 1)(z + 1) & 0 \\
0 & (x^2 - 1)(z - 1) & 0 \\
0 & 0 & (x^2 - 1)(y + 1) \\
0 & 0 & (x^2 - 1)(y - 1) \\
0 & 0 & (y^2 - 1)(x + 1) \\
0 & 0 & (y^2 - 1)(x - 1) \\
4(y^2 - 1)x(z + 1) & 0 & (x^2 - 1)(y^2 - 1) \\
4(y^2 - 1)x(z - 1) & 0 & (x^2 - 1)(y^2 - 1) \\
(y^2 - 1)y(z + 1) & 0 & 0 \\
(y^2 - 1)y(z - 1) & 0 & 0 \\
4(z^2 - 1)x(y + 1) & (x^2 - 1)(z^2 - 1) & 0 \\
4(z^2 - 1)x(y - 1) & (x^2 - 1)(z^2 - 1) & 0 \\
(z^2 - 1)(y + 1)z & 0 & 0 \\
(z^2 - 1)(y - 1)z & 0 & 0 \\
0 & 4(x^2 - 1)y(z + 1) & (x^2 - 1)(y^2 - 1) \\
0 & 4(x^2 - 1)y(z - 1) & (x^2 - 1)(y^2 - 1) \\
0 & (x^2 - 1)x(z + 1) & 0 \\
0 & (x^2 - 1)x(z - 1) & 0 \\
(y^2 - 1)(z^2 - 1) & 4(z^2 - 1)(x + 1)y & 0 \\
(y^2 - 1)(z^2 - 1) & 4(z^2 - 1)(x - 1)y & 0 \\
0 & (z^2 - 1)(x + 1)z & 0 \\
0 & (z^2 - 1)(x - 1)z & 0 \\
0 & (x^2 - 1)(z^2 - 1) & 4(x^2 - 1)(y + 1)z \\
0 & (x^2 - 1)(z^2 - 1) & 4(x^2 - 1)(y - 1)z \\
0 & 0 & (x^2 - 1)x(y + 1) \\
0 & 0 & (x^2 - 1)x(y - 1) \\
(y^2 - 1)(z^2 - 1) & 0 & 4(y^2 - 1)(x + 1)z \\
(y^2 - 1)(z^2 - 1) & 0 & 4(y^2 - 1)(x - 1)z \\
0 & 0 & (y^2 - 1)(x + 1)y \\
0 & 0 & (y^2 - 1)(x - 1)y \\
\end{array} 
\end{array} \]
\[
\begin{array}{cl}
\calS_1^-\Lambda^1(\square_3)
&
\begin{array}{lll}
dx & dy & dz \\
\hline
(y + 1)(z + 1) & 0 & 0 \\
(y + 1)(z - 1) & 0 & 0 \\
(y - 1)(z + 1) & 0 & 0 \\
(y - 1)(z - 1) & 0 & 0 \\
0 & (x + 1)(z + 1) & 0 \\
0 & (x - 1)(z + 1) & 0 \\
0 & (x + 1)(z - 1) & 0 \\
0 & (x - 1)(z - 1) & 0 \\
0 & 0 & (x + 1)(y + 1) \\
0 & 0 & (x + 1)(y - 1) \\
0 & 0 & (x - 1)(y + 1) \\
0 & 0 & (x - 1)(y - 1) \\
\end{array} 
\end{array} \]
\[
\begin{array}{cl}
\calS_2^-\Lambda^1(\square_3)
&
\begin{array}{lll}
dx & dy & dz \\
\hline
(y + 1)(z + 1) & 0 & 0 \\
(y + 1)(z - 1) & 0 & 0 \\
(y - 1)(z + 1) & 0 & 0 \\
(y - 1)(z - 1) & 0 & 0 \\
0 & (x + 1)(z + 1) & 0 \\
0 & (x - 1)(z + 1) & 0 \\
0 & (x + 1)(z - 1) & 0 \\
0 & (x - 1)(z - 1) & 0 \\
0 & 0 & (x + 1)(y + 1) \\
0 & 0 & (x + 1)(y - 1) \\
0 & 0 & (x - 1)(y + 1) \\
0 & 0 & (x - 1)(y - 1) \\
x(y + 1)(z + 1) & 0 & 0 \\
x(y + 1)(z - 1) & 0 & 0 \\
x(y - 1)(z + 1) & 0 & 0 \\
x(y - 1)(z - 1) & 0 & 0 \\
0 & (x + 1)y(z + 1) & 0 \\
0 & (x - 1)y(z + 1) & 0 \\
0 & (x + 1)y(z - 1) & 0 \\
0 & (x - 1)y(z - 1) & 0 \\
0 & 0 & (x + 1)(y + 1)z \\
0 & 0 & (x + 1)(y - 1)z \\
0 & 0 & (x - 1)(y + 1)z \\
0 & 0 & (x - 1)(y - 1)z \\
(y^2 - 1)(z + 1) & 0 & 0 \\
(y^2 - 1)(z - 1) & 0 & 0 \\
(z^2 - 1)(y + 1) & 0 & 0 \\
(z^2 - 1)(y - 1) & 0 & 0 \\
0 & (z^2 - 1)(x + 1) & 0 \\
0 & (z^2 - 1)(x - 1) & 0 \\
0 & (x^2 - 1)(z + 1) & 0 \\
0 & (x^2 - 1)(z - 1) & 0 \\
0 & 0 & (x^2 - 1)(y + 1) \\
0 & 0 & (x^2 - 1)(y - 1) \\
0 & 0 & (y^2 - 1)(x + 1) \\
0 & 0 & (y^2 - 1)(x - 1) \\
\end{array} 
\end{array} \]
\[
\begin{array}{cl}
\calS_3^-\Lambda^1(\square_3)
&
\begin{array}{lll}
dx & dy & dz \\
\hline
(y + 1)(z + 1) & 0 & 0 \\
(y + 1)(z - 1) & 0 & 0 \\
(y - 1)(z + 1) & 0 & 0 \\
(y - 1)(z - 1) & 0 & 0 \\
0 & (x + 1)(z + 1) & 0 \\
0 & (x - 1)(z + 1) & 0 \\
0 & (x + 1)(z - 1) & 0 \\
0 & (x - 1)(z - 1) & 0 \\
0 & 0 & (x + 1)(y + 1) \\
0 & 0 & (x + 1)(y - 1) \\
0 & 0 & (x - 1)(y + 1) \\
0 & 0 & (x - 1)(y - 1) \\
x(y + 1)(z + 1) & 0 & 0 \\
x(y + 1)(z - 1) & 0 & 0 \\
x(y - 1)(z + 1) & 0 & 0 \\
x(y - 1)(z - 1) & 0 & 0 \\
0 & (x + 1)y(z + 1) & 0 \\
0 & (x - 1)y(z + 1) & 0 \\
0 & (x + 1)y(z - 1) & 0 \\
0 & (x - 1)y(z - 1) & 0 \\
0 & 0 & (x + 1)(y + 1)z \\
0 & 0 & (x + 1)(y - 1)z \\
0 & 0 & (x - 1)(y + 1)z \\
0 & 0 & (x - 1)(y - 1)z \\
x^2(y + 1)(z + 1) & 0 & 0 \\
x^2(y + 1)(z - 1) & 0 & 0 \\
x^2(y - 1)(z + 1) & 0 & 0 \\
x^2(y - 1)(z - 1) & 0 & 0 \\
0 & (x + 1)y^2(z + 1) & 0 \\
0 & (x - 1)y^2(z + 1) & 0 \\
0 & (x + 1)y^2(z - 1) & 0 \\
0 & (x - 1)y^2(z - 1) & 0 \\
0 & 0 & (x + 1)(y + 1)z^2 \\
0 & 0 & (x + 1)(y - 1)z^2 \\
0 & 0 & (x - 1)(y + 1)z^2 \\
0 & 0 & (x - 1)(y - 1)z^2 \\
(y^2 - 1)(z + 1) & 0 & 0 \\
(y^2 - 1)(z - 1) & 0 & 0 \\
(z^2 - 1)(y + 1) & 0 & 0 \\
(z^2 - 1)(y - 1) & 0 & 0 \\
0 & (z^2 - 1)(x + 1) & 0 \\
0 & (z^2 - 1)(x - 1) & 0 \\
0 & (x^2 - 1)(z + 1) & 0 \\
0 & (x^2 - 1)(z - 1) & 0 \\
0 & 0 & (x^2 - 1)(y + 1) \\
0 & 0 & (x^2 - 1)(y - 1) \\
0 & 0 & (y^2 - 1)(x + 1) \\
0 & 0 & (y^2 - 1)(x - 1) \\
(y^2 - 1)y(z + 1) & 0 & 0 \\
(y^2 - 1)y(z - 1) & 0 & 0 \\
(z^2 - 1)(y + 1)z & 0 & 0 \\
(z^2 - 1)(y - 1)z & 0 & 0 \\
0 & (z^2 - 1)(x + 1)z & 0 \\
0 & (z^2 - 1)(x - 1)z & 0 \\
0 & (x^2 - 1)x(z + 1) & 0 \\
0 & (x^2 - 1)x(z - 1) & 0 \\
0 & 0 & (x^2 - 1)x(y + 1) \\
0 & 0 & (x^2 - 1)x(y - 1) \\
0 & 0 & (y^2 - 1)(x + 1)y \\
0 & 0 & (y^2 - 1)(x - 1)y \\
(y^2 - 1)x(z + 1) & -(x^2 - 1)y(z + 1) & 0 \\
(y^2 - 1)x(z - 1) & -(x^2 - 1)y(z - 1) & 0 \\
(z^2 - 1)x(y + 1) & 0 & -(x^2 - 1)(y + 1)z \\
(z^2 - 1)x(y - 1) & 0 & -(x^2 - 1)(y - 1)z \\
0 & (z^2 - 1)(x + 1)y & -(y^2 - 1)(x + 1)z \\
0 & (z^2 - 1)(x - 1)y & -(y^2 - 1)(x - 1)z \\
\end{array} 
\end{array}\]
\[
\begin{array}{cl}
\calS_1\Lambda^2(\square_3)
&
\begin{array}{lll}
dydz & dxdz & dxdy \\
\hline
x + 1 & 0 & 0 \\
x - 1 & 0 & 0 \\
0 & y + 1 & 0 \\
0 & y - 1 & 0 \\
0 & 0 & z + 1 \\
0 & 0 & z - 1 \\
2(x + 1)y & y^2 - 1 & 0 \\
2(x - 1)y & y^2 - 1 & 0 \\
2(x + 1)z & 0 & -z^2 + 1 \\
2(x - 1)z & 0 & -z^2 + 1 \\
x^2 - 1 & 2x(y + 1) & 0 \\
x^2 - 1 & 2x(y - 1) & 0 \\
0 & 2(y + 1)z & z^2 - 1 \\
0 & 2(y - 1)z & z^2 - 1 \\
-x^2 + 1 & 0 & 2x(z + 1) \\
-x^2 + 1 & 0 & 2x(z - 1) \\
0 & y^2 - 1 & 2y(z + 1) \\
0 & y^2 - 1 & 2y(z - 1) \\
\end{array} 
\end{array} \]
\[
\begin{array}{cl}
\calS_2\Lambda^2(\square_3)
&
\begin{array}{lll}
dydz & dxdz & dxdy \\
\hline
x + 1 & 0 & 0 \\
x - 1 & 0 & 0 \\
0 & y + 1 & 0 \\
0 & y - 1 & 0 \\
0 & 0 & z + 1 \\
0 & 0 & z - 1 \\
(x + 1)z & 0 & 0 \\
(x - 1)z & 0 & 0 \\
(x + 1)y & 0 & 0 \\
(x - 1)y & 0 & 0 \\
0 & (y + 1)z & 0 \\
0 & (y - 1)z & 0 \\
0 & x(y + 1) & 0 \\
0 & x(y - 1) & 0 \\
0 & 0 & y(z + 1) \\
0 & 0 & y(z - 1) \\
0 & 0 & x(z + 1) \\
0 & 0 & x(z - 1) \\
3(x + 1)y^2 & (y^2 - 1)y & 0 \\
3(x - 1)y^2 & (y^2 - 1)y & 0 \\
3(x + 1)z^2 & 0 & -(z^2 - 1)z \\
3(x - 1)z^2 & 0 & -(z^2 - 1)z \\
(x^2 - 1)x & 3x^2(y + 1) & 0 \\
(x^2 - 1)x & 3x^2(y - 1) & 0 \\
0 & 3(y + 1)z^2 & (z^2 - 1)z \\
0 & 3(y - 1)z^2 & (z^2 - 1)z \\
-(x^2 - 1)x & 0 & 3x^2(z + 1) \\
-(x^2 - 1)x & 0 & 3x^2(z - 1) \\
0 & (y^2 - 1)y & 3y^2(z + 1) \\
0 & (y^2 - 1)y & 3y^2(z - 1) \\
4(x + 1)yz & (y^2 - 1)z & -(z^2 - 1)y \\
4(x - 1)yz & (y^2 - 1)z & -(z^2 - 1)y \\
(x^2 - 1)z & 4x(y + 1)z & (z^2 - 1)x \\
(x^2 - 1)z & 4x(y - 1)z & (z^2 - 1)x \\
-(x^2 - 1)y & (y^2 - 1)x & 4xy(z + 1) \\
-(x^2 - 1)y & (y^2 - 1)x & 4xy(z - 1) \\
x^2 - 1 & 0 & 0 \\
0 & y^2 - 1 & 0 \\
0 & 0 & z^2 - 1 \\

\end{array} 
\end{array} \]
\[
\begin{array}{cl}
\calS_3\Lambda^2(\square_3)
&
\begin{array}{lll}
dydz & dxdz & dxdy \\
\hline
x + 1 & 0 & 0 \\
x - 1 & 0 & 0 \\
0 & y + 1 & 0 \\
0 & y - 1 & 0 \\
0 & 0 & z + 1 \\
0 & 0 & z - 1 \\
(x + 1)z & 0 & 0 \\
(x - 1)z & 0 & 0 \\
(x + 1)y & 0 & 0 \\
(x - 1)y & 0 & 0 \\
0 & (y + 1)z & 0 \\
0 & (y - 1)z & 0 \\
0 & x(y + 1) & 0 \\
0 & x(y - 1) & 0 \\
0 & 0 & y(z + 1) \\
0 & 0 & y(z - 1) \\
0 & 0 & x(z + 1) \\
0 & 0 & x(z - 1) \\
(x + 1)z^2 & 0 & 0 \\
(x - 1)z^2 & 0 & 0 \\
(x + 1)yz & 0 & 0 \\
(x - 1)yz & 0 & 0 \\
(x + 1)y^2 & 0 & 0 \\
(x - 1)y^2 & 0 & 0 \\
0 & (y + 1)z^2 & 0 \\
0 & (y - 1)z^2 & 0 \\
0 & x(y + 1)z & 0 \\
0 & x(y - 1)z & 0 \\
0 & x^2(y + 1) & 0 \\
0 & x^2(y - 1) & 0 \\
0 & 0 & y^2(z + 1) \\
0 & 0 & y^2(z - 1) \\
0 & 0 & xy(z + 1) \\
0 & 0 & xy(z - 1) \\
0 & 0 & x^2(z + 1) \\
0 & 0 & x^2(z - 1) \\
4(x + 1)y^3 & (y^2 - 1)y^2 & 0 \\
4(x - 1)y^3 & (y^2 - 1)y^2 & 0 \\
4(x + 1)z^3 & 0 & -(z^2 - 1)z^2 \\
4(x - 1)z^3 & 0 & -(z^2 - 1)z^2 \\
(x^2 - 1)x^2 & 4x^3(y + 1) & 0 \\
(x^2 - 1)x^2 & 4x^3(y - 1) & 0 \\
0 & 4(y + 1)z^3 & (z^2 - 1)z^2 \\
0 & 4(y - 1)z^3 & (z^2 - 1)z^2 \\
-(x^2 - 1)x^2 & 0 & 4x^3(z + 1) \\
-(x^2 - 1)x^2 & 0 & 4x^3(z - 1) \\
0 & (y^2 - 1)y^2 & 4y^3(z + 1) \\
0 & (y^2 - 1)y^2 & 4y^3(z - 1) \\
5(x + 1)yz^2 & (y^2 - 1)z^2 & -(z^2 - 1)yz \\
5(x - 1)yz^2 & (y^2 - 1)z^2 & -(z^2 - 1)yz \\
5(x + 1)y^2z & (y^2 - 1)yz & -(z^2 - 1)y^2 \\
5(x - 1)y^2z & (y^2 - 1)yz & -(z^2 - 1)y^2 \\
(x^2 - 1)xz & 5x^2(y + 1)z & (z^2 - 1)x^2 \\
(x^2 - 1)xz & 5x^2(y - 1)z & (z^2 - 1)x^2 \\
(x^2 - 1)z^2 & 5x(y + 1)z^2 & (z^2 - 1)xz \\
(x^2 - 1)z^2 & 5x(y - 1)z^2 & (z^2 - 1)xz \\
-(x^2 - 1)y^2 & (y^2 - 1)xy & 5xy^2(z + 1) \\
-(x^2 - 1)y^2 & (y^2 - 1)xy & 5xy^2(z - 1) \\
-(x^2 - 1)xy & (y^2 - 1)x^2 & 5x^2y(z + 1) \\
-(x^2 - 1)xy & (y^2 - 1)x^2 & 5x^2y(z - 1) \\
x^2 - 1 & 0 & 0 \\
0 & y^2 - 1 & 0 \\
0 & 0 & z^2 - 1 \\
(x^2 - 1)z & 0 & 0 \\
(x^2 - 1)y & 0 & 0 \\
(x^2 - 1)x & 0 & 0 \\
0 & (y^2 - 1)x & 0 \\
0 & (y^2 - 1)z & 0 \\
0 & (y^2 - 1)y & 0 \\
0 & 0 & (z^2 - 1)y \\
0 & 0 & (z^2 - 1)x \\
0 & 0 & (z^2 - 1)z \\
\end{array} 
\end{array} \]

\newpage
\[
\begin{array}{cl}
\calS_1^-\Lambda^2(\square_3)
&
\begin{array}{lll}
dydz & dxdz & dxdy \\
\hline
x + 1 & 0 & 0 \\
x - 1 & 0 & 0 \\
0 & y + 1 & 0 \\
0 & y - 1 & 0 \\
0 & 0 & z + 1 \\
0 & 0 & z - 1 \\
\end{array} 
\end{array} \]
\[
\begin{array}{cl}
\calS_2^-\Lambda^2(\square_3)
&
\begin{array}{lll}
dydz & dxdz & dxdy \\
\hline
x + 1 & 0 & 0 \\
x - 1 & 0 & 0 \\
0 & y + 1 & 0 \\
0 & y - 1 & 0 \\
0 & 0 & z + 1 \\
0 & 0 & z - 1 \\
(x + 1)z & 0 & 0 \\
(x - 1)z & 0 & 0 \\
(x + 1)y & 0 & 0 \\
(x - 1)y & 0 & 0 \\
0 & (y + 1)z & 0 \\
0 & (y - 1)z & 0 \\
0 & x(y + 1) & 0 \\
0 & x(y - 1) & 0 \\
0 & 0 & y(z + 1) \\
0 & 0 & y(z - 1) \\
0 & 0 & x(z + 1) \\
0 & 0 & x(z - 1) \\
x^2 - 1 & 0 & 0 \\
0 & y^2 - 1 & 0 \\
0 & 0 & z^2 - 1 \\
\end{array} 
\end{array} \]
\[
\begin{array}{cl}
\calS_3^-\Lambda^2(\square_3)
&
\begin{array}{lll}
dydz & dxdz & dxdy \\
\hline
x + 1 & 0 & 0 \\
x - 1 & 0 & 0 \\
0 & y + 1 & 0 \\
0 & y - 1 & 0 \\
0 & 0 & z + 1 \\
0 & 0 & z - 1 \\
(x + 1)z & 0 & 0 \\
(x - 1)z & 0 & 0 \\
(x + 1)y & 0 & 0 \\
(x - 1)y & 0 & 0 \\
0 & (y + 1)z & 0 \\
0 & (y - 1)z & 0 \\
0 & x(y + 1) & 0 \\
0 & x(y - 1) & 0 \\
0 & 0 & y(z + 1) \\
0 & 0 & y(z - 1) \\
0 & 0 & x(z + 1) \\
0 & 0 & x(z - 1) \\
(x + 1)z^2 & 0 & 0 \\
(x - 1)z^2 & 0 & 0 \\
(x + 1)yz & 0 & 0 \\
(x - 1)yz & 0 & 0 \\
(x + 1)y^2 & 0 & 0 \\
(x - 1)y^2 & 0 & 0 \\
0 & (y + 1)z^2 & 0 \\
0 & (y - 1)z^2 & 0 \\
0 & x(y + 1)z & 0 \\
0 & x(y - 1)z & 0 \\
0 & x^2(y + 1) & 0 \\
0 & x^2(y - 1) & 0 \\
0 & 0 & y^2(z + 1) \\
0 & 0 & y^2(z - 1) \\
0 & 0 & xy(z + 1) \\
0 & 0 & xy(z - 1) \\
0 & 0 & x^2(z + 1) \\
0 & 0 & x^2(z - 1) \\
x^2 - 1 & 0 & 0 \\
0 & y^2 - 1 & 0 \\
0 & 0 & z^2 - 1 \\
(x^2 - 1)x & 0 & 0 \\
0 & (y^2 - 1)y & 0 \\
0 & 0 & (z^2 - 1)z \\
(x^2 - 1)y & -(y^2 - 1)x & 0 \\
(x^2 - 1)z & 0 & (z^2 - 1)x \\
0 & (y^2 - 1)z & -(z^2 - 1)y \\
\end{array} 
\end{array} \]
\end{tiny}

\end{appendices}

\end{document}